\documentclass[12pt,a4paper]{article}
\usepackage[utf8]{inputenc}

\title{%
  Detecting changes in cross-sectional dependence in multivariate time series
}
\author{%
  Axel B\"ucher \\
  \small{Institut f\"ur Angewandte Mathematik} \\
  \small{Universit\"at Heidelberg} \\
  \small{Im Neuenheimer Feld 294, 69120 Heidelberg, Germany} \\
  \small{\texttt{axel.buecher@rub.de}}\\
  \and
  Ivan Kojadinovic, Tom Rohmer \\
  \small{Laboratoire de math\'ematiques et applications, UMR CNRS 5142} \\
  \small{Universit\'e de Pau et des Pays de l'Adour} \\
  \small{B.P. 1155, 64013 Pau Cedex, France} \\
  \small{\texttt{\{ivan.kojadinovic,tom.rohmer\}@univ-pau.fr}}
  \and
  Johan Segers\\ 
  \small{Institut de statistique, biostatistique et sciences actuarielles} \\
  \small{Université catholique de Louvain} \\
  \small{Voie du Roman Pays 20, B-1348 Louvain-la-Neuve, Belgium} \\
  \small{\texttt{johan.segers@uclouvain.be}}%
}

\usepackage{fullpage}
\usepackage{amssymb}
\usepackage{amsmath}
\usepackage{amsthm}
\usepackage[round]{natbib}
\usepackage{url}
\usepackage{graphicx}
\usepackage{multirow}
\usepackage{setspace}
\usepackage{enumerate}
\usepackage{paralist}
\usepackage{color}
\usepackage{booktabs,rotating}
\usepackage{bm}

\numberwithin{equation}{section}

\newcommand{\eps}{\varepsilon}
\newcommand{\N}{\mathbb{N}}
\newcommand{\R}{\mathbb{R}}

\newcommand{\dd}{\mathrm{d}}
\newcommand{\Cb}{\mathbb{C}}

\newcommand{\FF}{\mathcal{F}}

\newcommand{\D}{\mathbb{D}}
\newcommand{\B}{\mathbb{B}}
\newcommand{\Z}{\mathbb{Z}}
\newcommand{\Ex}{\operatorname{E}}
\newcommand{\expec}{\Ex}

\newcommand{\cov}{\operatorname{cov}}
\newcommand{\argmax}{\operatornamewithlimits{\arg\max}}
\newcommand{\1}{\mathbf{1}}
\newcommand{\ip}[1]{\lfloor #1 \rfloor}
\renewcommand{\vec}{\bm}
\newcommand{\pobs}[1]{\hat{\bm #1}}
\renewcommand{\Pr}{\mathrm{P}}
\newcommand{\p}{\overset{\Pr}{\to}}

\newcommand{\Bseq}{\B_n} % Ordinary sequential empirical process for U_i; Formerly \tilde{\B}_n
\newcommand{\Bseqcen}{\mathring{\B}_n} % auto-centered sequential empirical process for U_i; Formerly \check{\tilde{\B}}_n

\theoremstyle{plain}
\newtheorem{prop}{Proposition}[section]

\newtheorem{cor}[prop]{Corollary}

\newtheorem{cond}[prop]{Condition}

\begin{document}

\maketitle

\begin{abstract}
Classical and more recent tests for detecting distributional changes in multivariate time series often lack power against alternatives that involve changes in the cross-sectional dependence structure. To be able to detect such changes better, a test is introduced based on a recently studied variant of the sequential empirical copula process. In contrast to earlier attempts, ranks are computed with respect to relevant subsamples, with beneficial consequences for the sensitivity of the test. For the computation of p-values we propose a multiplier resampling scheme that takes the serial dependence into account. The large-sample theory for the test statistic and the resampling scheme is developed. The finite-sample performance of the procedure is assessed by Monte Carlo simulations. Two case studies involving time series of financial returns are presented as well.
\end{abstract}

%\noindent {\it Keywords:} change-point detection; empirical copula; multiplier central limit theorem; partial-sum process; ranks; strong mixing.

\section{Introduction}
% ====================

Given a sequence $\vec X_1,\dots,\vec X_n$ of $d$-dimensional observations, change-point detection aims at testing
\begin{equation}
\label{H0}
  H_0 : \,\exists \, F \text{ such that } 
  \vec X_1, \ldots, \vec X_n \text{ have c.d.f. } F
\end{equation}
against alternatives involving the nonconstancy of the c.d.f. Under $H_0$ and the assumption that $\vec X_1,\dots,\vec X_n$ have continuous marginal c.d.f.s $F_1,\dots,F_d$, we have from the work of \cite{Skl59} that the common multivariate c.d.f.\ $F$ can be written in a unique way as
$$
F(\vec x) = C \{ F_1(x_1),\dots,F_d(x_d) \}, \qquad \vec x \in \R^d,
$$
where the function $C:[0, 1]^d \to [0,1]$ is a \emph{copula} and can be regarded as capturing the dependence between the components of $\vec X_1,\dots,\vec X_n$. It follows that $H_0$ can be rewritten as $H_{0,m} \cap H_{0,c}$, where
\begin{align}
\label{H0m}
  H_{0,m} &: \,\exists \, F_1,\dots, F_d \text{ such that } 
  \vec X_1, \ldots, \vec X_n \text{ have marginal c.d.f.s } F_1, \dots, F_d,
\\
\label{H0c}
  H_{0,c} &: \,\exists \, C \text{ such that } 
  \vec X_1, \ldots, \vec X_n \text{ have copula } C.
\end{align}

Classical nonparametric tests for $H_0$ are based on sequential empirical processes; see e.g.\ \citet{Bai94}, \citet[Section~2.6]{CsoHor97} and \citet{Ino01}. For moderate sample sizes, however, such tests appear to have little power against alternative hypotheses that leave the margins unchanged but that involve a change in the copula, i.e., when $H_{0,m} \cap (\neg H_{0,c})$ holds. Empirical evidence of the latter fact can be found in \citet[Section~4]{HolKojQue13}.
% and will be additionally given in Section~\ref{simulations}. 
For that reason, nonparametric tests for change-point detection particularly sensitive to changes in the dependence structure are of practical interest.

Several tests designed to capture changes in cross-sectional dependence structure were proposed in the literature. Tests based on Kendall's tau were investigated by \cite{GomHor99} \citep[see also][]{GomHor02}, \cite{QueSaiFav13} and \cite{DehVogWenWie13}. Although these have good power when the copula changes in such a way that Kendall's tau changes as well, they are obviously useless when the copula changes but Kendall's tau does not change or only very little. Tests based on functionals of sequential empirical copula processes were considered in \citet{Rem10}, \cite{BucRup13}, \cite{vanWie13} and \cite{WieDehvanVog13}. However, the power of such tests is often disappointing; see Section~\ref{sec:simulations} for some numerical evidence.

It is our aim to construct a new test for $H_0$ that is more powerful than its predecessors against alternatives that involve a change in the copula. The test is based on sequential empirical copula processes as well, but the crucial difference lies in the computation of the ranks. Whereas in \citet{Rem10} and subsequent papers, ranks are always computed with respect to the full sample, we propose to compute the ranks with respect to the relevant subsamples; see Section~\ref{sec:test} for details. The intuition is that in this way, the copulas of those subsamples are estimated more accurately, so that differences between copulas of disjoint subsamples are detected more quickly. The phenomenon is akin to the one observed in \cite{GenSeg10} that the empirical copula, which is based on pseudo-observations, is often a better estimator of a copula than the empirical distribution function based on observations from the copula itself. For another illustration in the context of tail dependence functions, see \cite{
Buc13}.

The paper is organized as follows. The test statistic is presented in Section~\ref{sec:test}, and its asymptotic distribution under the null hypothesis is found in Section~\ref{sec:large-sample}. Next, Section~\ref{sec:resampling} contains a detailed description of the multiplier resampling scheme and its asymptotic validity under the null hypothesis. The results of a large-scale Monte Carlo simulation study are reported in Section~\ref{sec:simulations}, and two brief case studies are given in Section~\ref{sec:cases}. Section~\ref{sec:conclusion} concludes. Proofs and details regarding the simulation study are deferred to the Appendices.

In the rest of the paper, the arrow~`$\leadsto$' denotes weak convergence in the sense of Definition~1.3.3 in \cite{vanWel96}. Given a set $T$, let $\ell^\infty(T)$ denote the space of all bounded real-valued functions on $T$ equipped with the uniform metric.

\section{Test statistic}
% ======================
\label{sec:test}

We now describe our test statistic and highlight the difference with the one in \citet{Rem10} and \citet{BucRup13}. Let $\vec X_1, \ldots, \vec X_n$ be random vectors. For integers $1 \le k \le l \le n$, let $C_{k:l}$ be the empirical copula of the sample $\vec X_k, \ldots, \vec X_l$. Specifically,
\begin{equation}
\label{eq:Ckl}
  C_{k:l}(\vec{u}) 
  = \frac{1}{l-k+1} \sum_{i=k}^l \1(\pobs{U}_i^{k:l} \leq \vec{u}),
\end{equation}
for $\vec{u} \in [0, 1]^d$, where
\begin{equation}
\label{eq:pseudo}
  \pobs{U}_i^{k:l} 
  = \frac{1}{l-k+1} (R_{i1}^{k:l},\dots,R_{id}^{k:l}), \qquad 
  i \in \{ k, \dots, l \},
\end{equation}
with $R_{ij}^{k:l} = \sum_{t=k}^l \1( X_{tj} \le X_{ij} )$ the (maximal) rank of $X_{ij}$ among $X_{kj},\dots,X_{lj}$. (Because of serial dependence, there can be ties, even if the marginal distribution is continuous; think for instance of a moving maximum process.) An important point is that the ranks are computed within the subsample $\vec X_k, \ldots, \vec X_l$ and not within the whole sample $\vec X_1, \ldots, \vec X_n$. As we continue, we adopt the convention that $C_{k:l}=0$ if $k > l$.

Write $\Delta = \{ (s, t) \in [0, 1]^2: s \le t \}$. Let $\lambda_n(s, t) = (\ip{nt}-\ip{ns})/n$ for $(s, t) \in \Delta$. Our test statistic is based on the difference process, $\D_n$, defined by
\begin{equation}
\label{eq:Dn}
\D_n(s,\vec{u}) = \sqrt{n} \, \lambda_n(0,s) \, \lambda_n(s,1) \, \{ C_{1:\ip{ns}}(\vec{u}) - C_{\ip{ns}+1:n}(\vec{u}) \}
\end{equation}
for $(s,\vec{u}) \in [0,1]^{d+1}$. For every $s \in [0, 1]$, it gives a weighted difference between the empirical copulas at $\vec{u}$ of the first $\ip{ns}$ and the last $n - \ip{ns}$ points of the sample. Large absolute differences point in the direction of a change in the copula.

To aggregate over $\vec{u}$, we consider the Cram\'er--von Mises statistic
$$
  S_{n,k} = 
  \int_{[0, 1]^d} 
    \left\{ \D_n \left( k/n, \vec{u} \right) \right\}^2 \, 
  \dd C_{1:n}(\vec{u}), 
  \qquad k \in \{1,\dots,n-1\}.
$$
The test statistic for detecting changes in cross-sectional dependence is then
\begin{equation}
\label{eq:Sn}
  S_n 
  = \max_{1 \leq k \leq n-1} S_{n,k} 
  = \sup_{s \in [0,1]} 
  \int_{[0, 1]^d} 
    \{ \D_n ( s, \vec{u} ) \}^2 \,
  \dd C_{1:n}(\vec{u}).
\end{equation}
Other aggregating functions can be thought of too, leading for instance to Kolmogorov--Smirnov and Kuiper statistics. In numerical experiments, the resulting tests were found to be less powerful than the one based on the Cram\'er--von Mises statistic and hence are not considered further in this paper.

The null hypothesis of a constant distribution is rejected when $S_n$ is large. The p-values are determined by the null distribution of $S_n$, whose large-sample limit is derived in Section~\ref{sec:large-sample}. To estimate the p-values from the data, a multiplier bootstrap method is proposed in Section~\ref{sec:resampling}. 

Finally, if $H_0$ is rejected, there could be one or several abrupt or smooth changes in the joint distribution. Moreover, the change(s) could concern one or more marginal distributions, the copula, or both. In the case where there is just a single (abrupt) change-point $k^* \in \{1, \ldots, n-1\}$, one can for instance estimate it by
\begin{equation}
\label{eq:kstarn}
  k^\star_n = \argmax_{1 \leq k \leq n-1} S_{n,k}.
\end{equation}
% To establish the principle, consistency of this estimator in the case of a single change in the copula (and constant margins) is established in Appendix~\ref{app:changepoint}. 
We do not pursue the issue of single or multiple change-point estimation nor the diagnosis of the nature of the change-point. 

Our test statistic $S_n$ differs from the one considered in \citet[Section~5.2]{Rem10} and \citet[Section~3.2]{BucRup13} in the way the copulas of the subsamples $\vec X_k, \ldots, \vec X_l$ are estimated. Rather than the empirical copula $C_{k:l}$, these authors propose to use
\begin{equation}
\label{eq:Ckln}
  C_{k:l,n}(\vec{u}) 
  = \frac{1}{l-k+1} \sum_{i=k}^l \1(\pobs{U}_i^{1:n} \leq \vec{u}),
  \qquad \vec{u} \in [0, 1]^d,
\end{equation}
with the convention that $C_{k:l,n}=0$ if $k > l$. In comparison with $C_{k:l}$ in~\eqref{eq:Ckl}, the ranks for the subsample $\vec X_k, \ldots, \vec X_l$ are computed relative to the complete sample $\vec X_1, \ldots, \vec X_n$. The estimators $C_{k:l,n}$ yield the difference process
\begin{equation}
\label{DnR}
  \D_n^R(s,\vec{u}) 
  = \sqrt{n} \, \lambda_n(0,s) \, \lambda_n(s,1) \,
  \{ C_{1:\ip{ns},n}(\vec{u}) - C_{\ip{ns}+1:n,n}(\vec{u}) \}
\end{equation}
for $(s,\vec{u}) \in [0,1]^{d+1}$. The process $\D_n^R$ is to be compared with the process $\D_n$ in~\eqref{eq:Dn}. The difference lies in the use of $C_{k:l,n}$ rather than $C_{k:l}$. From the process $\D_n^R$, one obtains the test statistic
\begin{equation}
\label{SnR}
  S_n^R 
  = \sup_{s \in [0,1]} 
  \int_{[0, 1]^d} 
    \left\{ \D_n^R \left( s, \vec{u} \right) \right\}^2 \,
  \dd C_{1:n}(\vec{u}),
\end{equation}
which is the analogue of $S_n$ in~\eqref{eq:Sn}.

In the Monte Carlo simulation experiments (Section~\ref{sec:simulations}), we will see that $S_n$ is usually more powerful than $S_n^R$ for detecting changes in the cross-sectional copula. Intuitively, the reason is that the empirical copula $C_{k:l}$ in~\eqref{eq:Ckl} is often a better copula estimator than $C_{k:l,n}$ in~\eqref{eq:Ckln}. Note that $C_{k:l}$ is not only the empirical copula of $\vec X_k, \ldots, \vec X_l$, it is also equal to the empirical copula of $\pobs{U}_{k}^{1:n}, \ldots, \pobs{U}_{l}^{1:n}$, of which $C_{k:l,n}$ is the empirical distribution function.

In \cite{GenSeg10}, situations are identified where the empirical copula of an independent random sample drawn from a given bivariate copula has a lower asymptotic variance than the empirical distribution function of that sample. Of course, the situation here is different from the one in the cited paper: multivariate rather than bivariate, serial dependence rather than independence. But still, we suspect the same mechanisms to be active.

\section{Large-sample distribution}
% =================================
\label{sec:large-sample}

The asymptotic distribution under $H_0$ of our test statistic $S_n$ in~\eqref{eq:Sn} can be obtained by writing it as a functional of the two-sided sequential empirical copula process studied in \cite{BucKoj13}. Let $\vec X_1, \vec X_2, \ldots$ be a strictly stationary $d$-variate time series with stationary c.d.f.\ $F$ having continuous margins $F_1, \ldots, F_d$ and copula $C$. Recall $C_{k:l}$ in~\eqref{eq:Ckl} and $\pobs{U}_i^{k:l}$ in~\eqref{eq:pseudo}. The two-sided sequential empirical copula process, $\Cb_n$, is defined by
\begin{align}
\label{eq:secp2}
  \Cb_n(s, t, \vec{u}) 
  &= \sqrt{n} \, \lambda_n(s, t) \, \{ C_{\ip{ns}+1: \ip{nt}}(\vec{u})  - C(\vec{u}) \} \\
\label{eq:seqempcop}
  &= \frac{1}{\sqrt{n}} \sum_{i=\ip{ns}+1}^{\ip{nt}} \left\{  \1 ( \pobs{U}_i^{\ip{ns}+1:\ip{nt}} \leq \vec u) - C(\vec{u}) \right\},
\end{align}
for $(s, t, \vec{u}) \in \Delta \times [0, 1]^d$. The link of $\Cb_n$ to our test statistic $S_n$ in~\eqref{eq:Sn} is that, under $H_0$, the difference process $\D_n$ in~\eqref{eq:Dn} can be written as
\begin{equation}
\label{DnH0}
  \D_n(s,\vec{u}) 
  = \lambda_n(s,1) \, \Cb_n(0,s,\vec{u}) 
  - \lambda_n(0,s) \, \Cb_n(s,1,\vec{u}),
\end{equation}
for $(s,\vec{u}) \in [0,1]^{d+1}$.

Before focusing on the weak limit of the process $\D_n$ under $H_0$, let us briefly recall the notion of \emph{strongly mixing sequence}. For a sequence of $d$-dimensional random vectors $(\vec Y_i)_{i \in \Z}$, the $\sigma$-field generated by $(\vec Y_i)_{a \leq i \leq b}$, $a, b \in \Z \cup \{-\infty,+\infty \}$, is denoted by $\FF_a^b$. The strong mixing coefficients corresponding to the sequence $(\vec Y_i)_{i \in \Z}$ are defined by
$$
\alpha_r = \sup_{p \in \Z} \sup_{A \in \FF_{-\infty}^p,B\in \FF_{p+r}^{+\infty}} | P(A \cap B) - P(A) P(B) |
$$
for positive integer $r$. The sequence $(\vec Y_i)_{i \in \Z}$ is said to be \emph{strongly mixing} if $\alpha_r \to 0$ as $r \to \infty$. 

The weak limit of the two-sided empirical copula process $\Cb_n$ defined in~\eqref{eq:seqempcop} under strong mixing was established in \cite{BucKoj13} under the following two conditions:

\begin{cond}
\label{cond:noTies}
With probability one, there are no ties in each of the $d$ component series $X_{1j}, X_{2j}, \ldots$, where $j \in \{1, \ldots, d\}$.
\end{cond}

\begin{cond} 
\label{cond:pd}
For any $j \in \{1,\dots,d\}$, the partial derivatives $\dot C_j = \partial C/\partial u_j$ exist and are continuous on $V_j = \{ \vec{u} \in [0, 1]^d : u_j \in (0,1) \}$. 
\end{cond}

Condition~\ref{cond:noTies} was considered in \cite{BucKoj13} as continuity of the marginal distributions is \emph{not} sufficient to guarantee the absence of ties when the observations are serially dependent \cite[see e.g.][Example~4.2]{BucSeg13}. One of the contributions of this work is to show that it actually can be dispensed with. Condition~\ref{cond:pd} was proposed in \cite{Seg12} and is nonrestrictive in the sense that it is necessary for the candidate weak limit of $\Cb_n$ to exist pointwise and have continuous trajectories.

As we continue, for any $j \in \{1,\dots,d\}$, we define $\dot C_j$ to be zero on the set $\{ \vec{u} \in [0, 1]^d : u_j \in \{0,1\} \}$ \citep[see also][]{Seg12,BucVol13}. %It then follows that, under Condition~\ref{cond:pd}, $\dot C_j$ is defined on the whole of $[0, 1]^d$. 
Also, for any $j \in \{1,\dots,d\}$ and any $\vec{u} \in [0, 1]^d$, $\vec{u}^{(j)}$ is the vector of $[0, 1]^d$ defined by $u^{(j)}_i = u_j$ if $i = j$ and 1 otherwise. 

The weak convergence of the process $\Cb_n$ defined in~\eqref{eq:seqempcop} actually follows from that of the process
\begin{equation}
\label{eq:seqep}
  \Bseq(s, t, \vec{u}) 
  = \frac{1}{\sqrt{n}} \sum_{i=\ip{ns}+1}^{\ip{nt}} \{\1(\vec U_i \leq \vec{u}) - C(\vec{u}) \}, \qquad (s, t,\vec{u}) \in \Delta \times [0, 1]^d,
\end{equation}
where $\vec U_1,\dots,\vec U_n$ is the unobservable sample obtained from $\vec X_1, \dots, \vec X_n$ by the probability integral transforms $U_{ij} = F_j(X_{ij})$, $i \in \{1,\dots,n\}$, $j \in \{1,\dots,d\}$, and with the convention that $\Bseq(s, t, \cdot) = 0$ if $\ip{nt} - \ip{ns} = 0$. 

If $\vec U_1,\dots,\vec U_n$ is drawn from a strictly stationary sequence $(\vec U_i)_{i \in \Z}$ whose strong mixing coefficients satisfy $\alpha_r = O(r^{-a})$ with $a > 1$, we have from \cite{Buc13b} that $\Bseq(0,\cdot,\cdot)$ converges weakly in $\ell^\infty([0,1]^{d+1})$ to a tight centered Gaussian process $\Z_C$ with covariance function
$$
  \cov\{\Z_C(s,\vec{u}), \Z_C(t,\vec v)\} = \min(s,t) \sum_{k \in \Z} \cov\{ \1(\vec U_0 \leq \vec{u}), \1(\vec U_k \leq \vec v) \}.
$$
The latter is actually a consequence of Lemma~2 in \cite{Buc13b} stating that $\Bseq(0,\cdot,\cdot)$ is asymptotically uniformly equicontinuous in probability, which in turn implies that $\Z_C$ has continuous trajectories with probability one.
%known as a \emph{$C$-Kiefer-M\"uller} process. 
As a consequence of the continuous mapping theorem, $\Bseq \leadsto \B_C$ in $\ell^\infty(\Delta \times [0, 1]^d)$, where
\begin{equation}
\label{eq:BC}
  \B_C(s, t,\vec{u}) 
  = \Z_C(t,\vec{u}) - \Z_C(s,\vec{u}), 
  \qquad (s, t,\vec{u}) \in \Delta \times [0,1]^d.
\end{equation}

The following result is a consequence of Theorem~3.4 of \cite{BucKoj13} and the arguments used in the proof of Lemma~A.2 of \cite{BucSeg13}. Its proof is given in Appendix~\ref{proof:prop:weak_Cn_sm}.

\begin{prop}
\label{prop:weak_Cn_sm}
Let $\vec X_1,\dots,\vec X_n$ be drawn from a strictly stationary sequence $(\vec X_i)_{i \in \Z}$ with continuous margins and whose strong mixing coefficients satisfy $\alpha_r = O(r^{-a})$, $a > 1$. Then, provided Condition~\ref{cond:pd} holds,
\begin{equation}
\label{eq:asymequivCn}
  \sup_{(s, t,\vec{u})\in \Delta \times [0, 1]^d} 
  \left| \Cb_n(s, t, \vec{u}) - \tilde \Cb_n(s, t, \vec{u}) \right| 
  \p 0,
\end{equation}
where
\begin{equation}
\label{eq:tildeCn}
\tilde \Cb_n(s, t, \vec{u}) 
= \Bseq(s, t, \vec{u}) 
- \sum_{j=1}^d \dot C_j(\vec{u}) \, \Bseq(s, t, \vec{u}^{(j)}).
\end{equation}
Consequently, $\Cb_n \leadsto \Cb_C$ in $\ell^\infty(\Delta \times [0, 1]^d)$, where, for $(s, t, \vec{u}) \in \Delta \times [0, 1]^d$,
\begin{equation}
\label{eq:ec}
\Cb_C(s, t, \vec{u}) 
= \B_C(s, t, \vec{u}) 
- \sum_{j=1}^d \dot C_j(\vec{u}) \, \B_C(s, t, \vec{u}^{(j)}).
\end{equation}
\end{prop}

In view of~\eqref{DnH0}, the weak limit of $\D_n$ under $H_0$ is a mere corollary of Proposition~\ref{prop:weak_Cn_sm} and the continuous mapping theorem.

\begin{cor}
Under the conditions of Proposition~\ref{prop:weak_Cn_sm}, $\D_n \leadsto \D_C$ in $\ell^\infty([0,1]^{d+1})$, where, for any $(s,\vec{u}) \in [0,1]^{d+1}$, 
\begin{equation}
  \label{Dc}
  \D_C(s,\vec{u}) =  \Cb_C(0,s,\vec{u}) - s \, \Cb_C(0,1,\vec{u}),
\end{equation}
with $\Cb_C$ defined in~\eqref{eq:ec}. As a consequence,
\begin{equation}
\label{eq:S}
  S_n \leadsto S = 
  \sup_{s \in [0,1]} 
  \int_{[0, 1]^d} 
    \{ \D_C(s,\vec{u}) \}^2 \,
  \dd C(\vec{u}).
\end{equation}
\end{cor}

The covariance function of $\D_C$ can be expressed in terms of the one of $\Cb_C$ by
\begin{equation*}
% \label{covDC}
  \cov \{ \D_C(s, \vec{u}), \D_C(t, \vec v) \} 
  = \{ \min(s, t) - st \} \, \cov \{ \Cb_C(0, 1, \vec{u}), \Cb_C(0, 1, \vec v) \}.
\end{equation*}

\section{Resampling}
% ==================
\label{sec:resampling}

In order to compute p-values for $S_n$ based on~\eqref{eq:S}, we propose to use resampling methods. Tracing back the definition of $S_n$ via $\D_n$ to $\Cb_n$ in~\eqref{eq:secp2} and using the approximation via $\tilde{\Cb}_n$ in~\eqref{eq:tildeCn}, we find that it suffices to construct a resampling scheme for $\Bseq$ defined in~\eqref{eq:seqep} and to estimate the first-order partial derivatives, $\dot{C}_j$, of $C$.

% \subsection{Estimation of partial derivatives}
% % --------------------------------------------

\subsection{Multiplier sequences}
% -------------------------------

In the case of i.i.d.\ observations, \cite{Sca05} proposed to use a \emph{multiplier} approach in the spirit of \citet[Chapter 2.9]{vanWel96} to resample $\Bseq$. When the first-order partial derivatives of $C$ are estimated by finite-differencing as in \cite{RemSca09}, the resulting resampling scheme for $\Cb_n$ is frequently referred to as a \emph{multiplier bootstrap}. In a nonsequential setting based on independent observations, \cite{BucDet10} compared the finite-sample behavior of the various resampling techniques proposed in the literature and concluded that the multiplier bootstrap of \cite{RemSca09} has, overall, the best finite-sample properties. This technique was revisited theoretically by \cite{Seg12} who showed its asymptotic validity under Condition~\ref{cond:pd}. A sequential generalization of the latter result will be stated later in this section. In the case of independent observations, the multiplier bootstrap is based on \emph{i.i.d.\ multiplier sequences}. We say that a sequence of 
random variables $(\xi_{i,n})_{i \in \Z}$ is an \emph{i.i.d.\ multiplier sequence} if:
\begin{enumerate}[({M}0)]
\item %\label{item:iid} 
$(\xi_{i,n})_{i \in \Z}$ is i.i.d., independent of $\vec{X}_1, \ldots, \vec{X}_n$, with distribution not changing with~$n$, having mean 0, variance 1, and being such that $\int_0^\infty \{ \Pr(|\xi_{0,n}| > x) \}^{1/2} \dd x < \infty$.
\end{enumerate}

Starting from the seminal work of \citet[Section~3.3]{Buh93}, \cite{BucRup13} and \cite{BucKoj13} have studied a \emph{dependent} multiplier bootstrap for $\Cb_n$ which extends the multiplier bootstrap of \cite{RemSca09} to the sequential and strongly mixing setting. The key idea in \cite{Buh93} is to replace i.i.d.\ multipliers by suitably serially dependent multipliers that will capture the serial dependence in the data. In the rest of the paper, we say that a sequence of random variables $(\xi_{i,n})_{i \in \Z}$ is a \emph{dependent multiplier sequence} if:
\begin{enumerate}[({M}1)]
\item %\label{item:moments} 
The sequence $(\xi_{i,n})_{i \in \Z}$ is strictly stationary with $\Ex(\xi_{0,n}) = 0$, $\Ex(\xi_{0,n}^2) = 1$ and $\sup_{n \geq 1} \Ex(|\xi_{0,n}|^\nu) < \infty$ for all $\nu \geq 1$, and is independent of the available sample $\vec X_1,\dots,\vec X_n$.
\item %\label{item:ln} 
There exists a sequence $\ell_n \to \infty$ of strictly positive constants such that $\ell_n = o(n)$ and the sequence $(\xi_{i,n})_{i \in \Z}$ is $\ell_n$-dependent, i.e., $\xi_{i,n}$ is independent of $\xi_{i+h,n}$ for all $h > \ell_n$ and $i \in \N$. 
\item %\label{item:varphi} 
There exists a function $\varphi:\R \to [0,1]$, symmetric around 0, continuous at $0$, satisfying $\varphi(0)=1$ and $\varphi(x)=0$ for all $|x| > 1$ such that $\Ex(\xi_{0,n} \xi_{h,n}) = \varphi(h/\ell_n)$ for all $h \in \Z$.
\end{enumerate}

Ways to generate dependent multiplier sequences are mentioned in Section~\ref{sec:simulations} and Appendix~\ref{strongmixing}. 

% \js{Say that $\ell_n$ will be estimated from the data -- see simulations and cite \cite{BucKoj13}.}

%Roughly speaking, starting from a large number of independent multiplier sequences, the multiplier bootstrap can be used to compute the same number of ``approximate independent copies'' of $\Cb_n$. The conditions imposed on the multiplier sequences, either~(M0) or~(M1)--(M3), will depend on whether the data are serially independent or just strongly mixing. We detail this construction below.

\subsection{Computing p-values via resampling}
% --------------------------------------------

Let $M$ be a large integer and let $(\xi_{i,n}^{(1)})_{i \in \Z},\dots,(\xi_{i,n}^{(M)})_{i \in \Z}$ be $M$ independent copies of the same multiplier sequence. We will define two multiplier resampling schemes for the process $\Bseq$ in~\eqref{eq:seqep}. These will lead to two resampling schemes for the test statistic $S_n$ in~\eqref{eq:Sn}, on the basis of which approximate p-values can be computed.

Recall that $C_{k:l}$ in~\eqref{eq:Ckl} is the empirical copula of $\vec{X}_k, \ldots, \vec{X}_l$, which is the empirical distribution of the vectors of rescaled ranks $\hat{\vec{U}}_i^{k:l}$ in~\eqref{eq:pseudo}. For any $m \in \{1,\dots,M\}$ and $(s, t,\vec{u}) \in \Delta \times [0, 1]^d$, let
\begin{equation}
\label{eq:hatBnm}
\hat{\B}_n^{(m)}(s, t,\vec{u}) = \frac{1}{\sqrt{n}} \sum_{i=\ip{ns}+1}^{\ip{nt}} \xi_{i,n}^{(m)} \{ \1 ( \pobs{U}_i^{1:n} \leq \vec{u} ) - C_{1:n}(\vec{u}) \},
\end{equation}
and
\begin{align}
\nonumber
\check{\B}_n^{(m)}(s, t,\vec{u}) &= \frac{1}{\sqrt{n}} \sum_{i=\ip{ns}+1}^{\ip{nt}} \xi_{i,n}^{(m)} \{ \1 ( \pobs{U}_i^{\ip{ns}+1:\ip{nt}} \leq \vec{u} ) - C_{\ip{ns}+1:\ip{nt}}(\vec{u}) \} \\ 
\label{eq:checkBnm}
&= \frac{1}{\sqrt{n}} \sum_{i=\ip{ns}+1}^{\ip{nt}} ( \xi_{i,n}^{(m)} - \bar  \xi_{\ip{ns}+1:\ip{nt}}^{(m)} ) \1 ( \pobs{U}_i^{\ip{ns}+1:\ip{nt}} \leq \vec{u} ),
\end{align}
where $\bar \xi_{k:l}^{(m)}$ is the arithmetic mean of $\xi_{i,n}^{(m)}$ for $i \in \{ k, \ldots, l \}$.
% where $\bar \xi_{\ip{ns}+1:\ip{nt}}^{(m)} = (\ip{nt} - \ip{ns})^{-1} \sum_{i=\ip{ns}+1}^{\ip{nt}} \xi_{i,n}^{(m)}$. 
By convention, the sums are zero if $\ip{ns} = \ip{nt}$. Note that the ranks are computed relative to the complete sample $\vec{X}_1, \ldots, \vec{X}_n$ for $\hat{\B}_n^{(m)}(s, t, \cdot \, )$, whereas they are computed relative to the subsample $\vec{X}_{\ip{ns}+1}, \ldots, \vec{X}_{\ip{nt}}$ for $\check{\B}_n^{(m)}(s, t, \cdot \, )$.

In order to get to resampling versions of $\tilde{\Cb}_n$ in~\eqref{eq:tildeCn}, we need estimators of the first-order partial derivatives of $C$. A simple estimator based on $\vec{X}_k, \ldots, \vec{X}_l$ consists of finite differencing at a bandwidth of $h \equiv h(k, l) = \min\{ (l-k+1)^{-1/2}, 1/2 \}$. Varying slightly upon the definition in \cite{RemSca09} and following \citet[Section~3]{KojSegYan11}, we put
$$
\dot{C}_{j,k:l}(\vec{u})
=
\frac{C_{k:l}( \vec{u} + h \vec{e}_j ) - C_{k:l}( \vec{u} - h \vec{e}_j )}{\min(u_j+h,1) - \max(u_j-h, 0)}
$$
for $\vec{u} \in [0, 1]^d$, where $\vec{e}_j$ is the $j$th canonical unit vector in $\R^d$. Note that if $h \le u_j \le 1-h$, the denominator is just $2h$. The more general form of the denominator corrects for boundary effects ($u_j$ close to $0$ or $1$). Proceeding for instance as in \citet[proof of Proposition 2]{KojSegYan11}, we find that the previous estimator is uniformly bounded.

The resampling versions $\hat{\B}_n^{(m)}$ and $\check{\B}_n^{(m)}$ of $\Bseq$ then lead to the following resampling versions for $\tilde{\Cb}_n$: for $(s, t,\vec{u}) \in \Delta \times [0, 1]^d$,
\begin{align}
\nonumber%\label{eq:hatCbnm}
  \hat{\Cb}_n^{(m)}(s, t, \vec{u}) 
  &= \hat{\B}_n^{(m)}(s, t, \vec{u}) - 
  \sum_{j=1}^d \dot C_{j,1:n}(\vec{u}) \,
  \hat{\B}_n^{(m)}(s, t, \vec{u}^{(j)}), \\
\label{eq:checkCbnm}
  \check{\Cb}_n^{(m)}(s, t, \vec{u}) 
  &= \check{\B}_n^{(m)}(s, t, \vec{u}) - 
  \sum_{j=1}^d \dot C_{j,\ip{ns}+1:\ip{nt}}(\vec{u}) \, 
  \check{\B}_n^{(m)}(s, t, \vec{u}^{(j)}).
\end{align}
Recall that $\lambda_n(s, t) = (\ip{nt} - \ip{ns})/n$. The difference process $\D_n$ is to be resampled by one of the following two methods:
\begin{align*}
\nonumber
  \hat{\D}_n^{(m)}(s,\vec{u}) 
  &= \lambda_n(s,1) \, \hat{\Cb}_n^{(m)}(0, s, \vec{u}) 
  - \lambda_n(0,s) \, \hat{\Cb}_n^{(m)}(s, 1, \vec{u}) \\ 
%\label{eq:hatDnm}
  &= \hat{\Cb}_n^{(m)}(0,s,\vec{u}) 
  - \lambda_n(0,s) \, \hat{\Cb}_n^{(m)}(0, 1, \vec{u}), \\[1em]
%\label{eq:checkDnm}
  \check{\D}_n^{(m)}(s,\vec{u}) 
  &= \lambda_n(s,1) \, \check{\Cb}_n^{(m)}(0, s, \vec{u}) 
  - \lambda_n(0,s) \, \check{\Cb}_n^{(m)}(s,1, \vec{u}).
\end{align*}
For resampling the test statistic, one has the choice between
\begin{align}
\label{eq:Snm:hat}
\hat{S}_n^{(m)} 
&= \sup_{s \in [0,1]} 
\int_{[0, 1]^d} 
  \{ \hat{\D}_n^{(m)} (s, \vec{u}) \}^2 \,
\dd C_{1:n}(\vec{u}), \\
\label{eq:Snm:check}
\check{S}_n^{(m)} 
&= \sup_{s \in [0,1]} 
\int_{[0, 1]^d} 
  \{ \check{\D}_n^{(m)} (s, \vec{u}) \}^2 \,
\dd C_{1:n}(\vec{u}).
\end{align}
Finally, approximate p-values of the observed test statistic $S_n$ can be computed via either
\begin{equation}
\label{eq:pval}
  \frac{1}{M} \sum_{m=1}^M 
  \1 \left( \hat{S}_n^{(m)} \geq S_n \right) 
  \qquad \text{or} \qquad 
  \frac{1}{M} \sum_{m=1}^M 
  \1 \left( \check{S}_n^{(m)} \geq S_n \right).
\end{equation}
The null hypothesis is rejected if the estimated p-value is smaller than the desired significance level.

% \js{Briefly explain how to compute p-values for $S_n^R$.}

By comparison, note that for the test statistic $S_n^R$ in~\eqref{SnR} based on the process $\D_n^R$ in~\eqref{DnR}, an approximate p-value can be computed using the multiplier processes
\begin{equation}
\label{eq:DnRm}
  \D_n^{R,(m)}(s,\vec{u}) 
  = \hat{\B}_n^{(m)}(0,s,\vec{u}) - \lambda_n(0, s) \, \hat{\B}_n^{(m)}(0,1,\vec{u}), 
%   \qquad (s, \vec{u}) \in [0,1]^{d+1}, \, m \in \{1, \dots, M\},
\end{equation}
where $\hat{\B}_n^{(m)}$ is defined in~\eqref{eq:hatBnm}; see also \citet[Section~5.2]{Rem10} and \citet[Section~3.2]{BucRup13}.

\subsection{Asymptotic validity of the resampling scheme}
% -------------------------------------------------------

We establish the asymptotic validity of the multiplier resampling schemes described above under the null hypothesis. First, we need to impose conditions on the data generating process $\vec{X}_1, \ldots, \vec{X}_n$ and the multiplier sequences $(\xi_{i,n}^{(m)})_{i \in \Z}$ for $m \in \{1, \ldots, M\}$.

\begin{cond}
\label{cond:DGP}
One of the following two conditions holds:
\begin{enumerate}[\bf (i)]
\item The random vectors $\vec X_1,\dots,\vec X_n$ are i.i.d.\ and $(\xi_{i,n}^{(1)})_{i \in \Z},\dots,(\xi_{i,n}^{(M)})_{i \in \Z}$ are independent copies of a multiplier sequence satisfying~(M0).
\item The random vectors $\vec X_1,\dots,\vec X_n$ are drawn from a strictly stationary sequence $(\vec X_i)_{i \in \Z}$ whose strong mixing coefficients satisfy $\alpha_r = O(r^{-a})$ for some $a > 3 + 3d/2$, and $(\xi_{i,n}^{(1)})_{i \in \Z},\dots,(\xi_{i,n}^{(M)})_{i \in \Z}$ are independent copies of a dependent multiplier sequence satisfying~(M1)--(M3) with $\ell_n = O(n^{1/2 - \gamma})$ for some $0 < \gamma < 1/2$. 
\end{enumerate}
In both cases, the stationary distribution of $\vec{X}_i$ has continuous margins and a copula $C$ satisfying Condition~\ref{cond:pd}.
\end{cond}

If the random vectors $\vec X_1, \ldots, \vec X_n$ are i.i.d., they can also be considered to be drawn from a strongly mixing, strictly stationary sequence. Hence, for the multiplier sequences $(\xi_{i,n}^{(m)})_{i \in \Z}$, one could either assume~(M0) or~(M1)--(M3): both should work. However, as discussed in \citet[Section~2]{BucKoj13}, the use of dependent multipliers in the case of independent observations is likely to result in an efficiency loss. This is illustrated in the Monte Carlo simulations reported in \citet[Section~3]{BucRup13} and carried out for the test based on the statistic $S_n^R$ defined in~\eqref{SnR} which is resampled using multiplier processes asymptotically equivalent to those given in~\eqref{eq:DnRm}: the use of dependent multipliers in the case of serially independent data usually results in a loss of power and in a slightly more conservative test. Thus, in finite samples, if there is no evidence against serial independence, it appears more sensible to work under~(M0). 

% Second, as the resampling scheme involves the estimation of the partial derivatives of $C$, we need to assume that these estimators are consistent. We will need two conditions, the second one being stronger than the first one.
% 
% \begin{cond}
% \label{cond:estpd} 
% There exists a constant $\kappa > 0$ such that $|\dot C_{j,1:n}(\vec{u})| \leq \kappa$ for all $j \in \{1,\dots,d\}$, $n \geq 1$ and $\vec{u} \in [0, 1]^d$, and, for any $\delta \in (0,1/2)$ and $j \in \{1,\dots,d\}$, 
% \begin{equation}
% \label{eq:Cj1nconsistent}
%   \sup_{\substack{\vec{u} \in [0, 1]^d \\ u_j \in [\delta,1-\delta]}}
%   |\dot C_{j,1:n}(\vec{u}) - \dot C_j(\vec{u})| 
%   \p 0. 
% \end{equation}
% \end{cond}
% 
% \begin{cond}
% \label{cond:estpd2} 
% Idem as Condition~\ref{cond:estpd}, but with convergence in probability in~\eqref{eq:Cj1nconsistent} replaced by almost sure convergence. \js{We might need to enforce the condition in order to complete the proof of Proposition~\ref{prop:check}.}
% \end{cond}
% 
% \js{Add a sentence that, at least in the i.i.d.\ case, we can prove that the finite-differencing estimator satisfies the previous condition. By the way, can't we prove it in the mixing case using convergence of the sequential empirical copula process?}

We can now state the asymptotic distributions of the multiplier resampling schemes under the null hypothesis of a constant distribution. We provide two propositions, one for the resampling scheme based on $\hat{\B}_n^{(m)}$ in~\eqref{eq:hatBnm} and another one for the scheme based on $\check{\B}_n^{(m)}$ in~\eqref{eq:checkBnm}.

\begin{prop}
\label{prop:hat}
If Condition~\ref{cond:DGP} holds, then
\[
  \left(\Cb_n, \hat{\Cb}_n^{(1)}, \dots, \hat{\Cb}_n^{(M)} \right) 
  \leadsto 
  \left(\Cb_C, \Cb_C^{(1)}, \dots, \Cb_C^{(M)} \right)
\]
in $\{\ell^\infty(\Delta \times [0, 1]^d)\}^{M+1}$, where $\Cb_C$ is defined in~\eqref{eq:ec}, and $\Cb_C^{(1)},\dots,\Cb_C^{(M)}$ are independent copies of $\Cb_C$. As a consequence, also
\[
  \left( \D_n, \hat{\D}_n^{(1)}, \dots, \hat{\D}_n^{(M)} \right) 
  \leadsto 
  \left( \D_C, \D_C^{(1)}, \dots, \D_C^{(M)} \right),
\]
in $\{ \ell^\infty([0,1]^{(d+1)}) \}^{M+1}$, where $\D_C$ is defined in~\eqref{Dc} and $\D_C^{(1)}, \dots,  \D_C^{(M)}$ are independent copies of $\D_C$. Finally,
\[
  \left( S_n, \hat{S}_n^{(1)}, \dots, \hat{S}_n^{(M)} \right) 
  \leadsto 
  \left( S,S^{(1)}, \dots, S^{(M)} \right)
\]
where $S$ is defined in~\eqref{eq:S} and $S^{(1)},\dots,S^{(M)}$ are independent copies of $S$.
\end{prop}

Under Condition~\ref{cond:DGP}(i), the above result can be easily proved by starting from Theorem~1 of \cite{HolKojQue13} and adapting the arguments used in \citet[proof of Proposition 3.2]{Seg12}. Under Conditions~\ref{cond:DGP}(ii) and~\ref{cond:noTies}, the result was obtained in \citet[Proposition 4.2]{BucKoj13}. The additional arguments allowing to avoid Condition~\ref{cond:noTies} will be given in the proof of the next result. 

\begin{prop}
\label{prop:check}
If Condition~\ref{cond:DGP} holds, then the conclusions of Proposition~\ref{prop:hat} also hold with $\hat{\Cb}_n^{(m)}$ replaced by $\check{\Cb}_n^{(m)}$, $\hat{\D}_n^{(m)}$ replaced by $\check{\D}_n^{(m)}$, and $\hat{S}_n^{(m)}$ replaced by $\check{S}_n^{(m)}$.
\end{prop}

The proof of Proposition~\ref{prop:check} is somewhat involved and is given in detail in Appendix~\ref{proof:prop:check}.

Combining the last claims of Propositions~\ref{prop:hat} and~\ref{prop:check} with Proposition F.1 in \cite{BucKoj13}, we obtain that a test based on $S_n$ whose p-value is computed using one of the two approaches in~\eqref{eq:pval} will hold its level asymptotically as $n \to \infty$ followed by $M \to \infty$.

\section{Simulation study}
% ========================
\label{sec:simulations}

Large-scale Monte Carlo experiments were carried out in order to study the finite-sample performance of the derived tests for detecting changes in cross-sectional dependence. The main questions addressed by the study are the following:
\begin{enumerate}[(i)]
\item How well do the tests hold their size under the null hypothesis $H_0$ in~\eqref{H0} of no change?
\item What is the power of the tests against the alternative $H_{1,c}$ of a single change in cross-sectional dependence at constant margins? Specifically, the alternative hypothesis is $H_{1,c} \cap H_{0,m}$ with $H_{0,m}$ in~\eqref{H0m} and $H_{1,c}$ defined by
\begin{align}
\nonumber
  H_{1,c} :\, &
  \exists \text{ distinct } C_1 \text { and } C_2\text{, and } 
  k^\star \in \{1, \ldots, n-1\} \text{ such that }\\
\label{H1c}
  &\vec X_1, \ldots, \vec X_{k^\star} \text{ have copula } C_1 \text{ and } \vec X_{k^\star+1}, \ldots, \vec X_n \text{ have copula } C_2.
\end{align}
\item What happens if the change in distribution is only due to a change in the margins, the copula remaining constant? Specifically, the alternative hypothesis is $H_{1,m} \cap H_{0,c}$ with $H_{0,c}$ given in~\eqref{H0c} and $H_{1,m}$ defined by
\begin{align}
\nonumber
H_{1,m} :\, &\exists \mbox{ distinct } F_{1,1}, F_{1,2} \mbox { as well as  } F_2,\dots,F_d\mbox{ and } k_1^\star \in \{1, \ldots, n-1\} \\
\nonumber
&\mbox{such that } \vec X_1, \ldots, \vec X_{k_1^\star} \mbox{ have marginal c.d.f.s } F_{1,1},F_2,\dots,F_d \\
\label{H1m}
&\mbox{and } \vec X_{k_1^\star+1}, \ldots, \vec X_n \mbox{ have marginal c.d.f.s } F_{1,2},F_2,\dots,F_d.
\end{align}
%Up to what extent can the tests be considered as a test for detecting changes in the copula? 
\end{enumerate}

In addition to the three questions above, many others can be formulated, involving other alternative hypotheses for instance. The problem is complex and there are countless ways of combining factors in the experimental design. In our study, the settings were chosen to represent a wide and hopefully representative variety of situations, in function of the three questions above. The main factors of our experiments are summarized below:
\begin{compactitem}
\item Test statistics: 
\begin{compactitem}
\item Our statistic $S_n$ in~\eqref{eq:Sn} with p-values computed via resampling using $\hat{S}_n^{(m)}$ or $\check{S}_n^{(m)}$ in~\eqref{eq:Snm:hat} and~\eqref{eq:Snm:check}, respectively. As we continue, we shall simply talk about the test based on $\hat S_n$ or $\check S_n$, respectively, to distinguish between these two situations.
\item The statistic $S_n^R$ in~\eqref{SnR} of \cite{BucRup13}, with p-values computed according to the resampling method for $\D_n^R$ in~\eqref{eq:DnRm}.
\end{compactitem}
\item Sample size: $n \in \{ 50, 100, 200 \}$.
\item Number of samples per setting: $1\,000$.
\item Cross-sectional dimension: $d \in \{ 2, 3 \}$.
\item Significance level: $\alpha = 5 \%$.
\item Serial dependence: The data were generated either as being serially independent or via two time-series models, an autoregressive process and a multivariate version of the exponential autoregressive model considered in \cite{AueTjo90} and \citet[Section~3.3]{PapPol01}. Independent standard normals were used as multipliers for independent observations, while for the serially dependent datasets, the dependent multiplier sequences were generated from initial independent standard normal sequences using the ``moving average approach'' proposed initially in \citet{Buh93} and revisited in some detail in \citet[Section~6.1]{BucKoj13}. The value of the bandwidth parameter $\ell_n$ defined in Condition~(M2) was chosen automatically using the approach described in \citet[Section~5]{BucKoj13}. See Appendix~\ref{app:tables} for details.
\item Margins: in all but one setting, the margins were kept constant, i.e., $H_{0,m}$ in \eqref{H0m} was assumed. In one case (see Table~\ref{H0break}), a break as in $H_{1,m}$ in \eqref{H1m} was assumed, the marginal distribution of the first component changing from the $N(0, 1)$ to the $N(\mu, 1)$ distribution.
\item Copulas: Clayton, Gumbel--Hougaard, Normal, Frank, with positive or negative (insofar possible) association, as well as asymmetric versions obtained via Khoudraji's device \citep{Kho95,GenGhoRiv98,Lie08}.
\item Alternative hypotheses involving a single change-point occurring at time $k^\star = \ip{nt}$ with $t \in \{0.1,0.25, 0.5, 0.75\}$:
\begin{compactitem}
\item $H_{0,m} \cap H_{1,c}$ with a change of the parameter within a copula family.
\item $H_{0,m} \cap H_{1,c}$ with a change of the copula family at constant Kendall's tau.
\item $H_{1,m} \cap H_{0,c}$, i.e., a change of one of the \emph{margins} rather than of the copula.
\item For the serially dependent case, a change in the copula of the innovations, leading to a gradual change of the copula of the marginal distributions of the observables.
\end{compactitem}
\end{compactitem}

The experiments were carried out in the \textsf{R} statistical system \citep{Rsystem} using the \texttt{copula} package \citep{copula}. To allow us to reuse previously written code, the rescaled ranks in~\eqref{eq:pseudo} were computed by dividing the ranks by $l-k+2$ instead of $l-k+1$. Because~\eqref{eq:Snm:hat} only involves rescaled ranks computed from the entire sample, the test based on $\hat S_n$ can be implemented to be substantially faster than the one based on $\check S_n$ for larger sample sizes. The corresponding routines are available in the \textsf{R} package \texttt{npcp} \citep{npcp}.

For the sake of brevity, only a representative subset of the results is reported here. Specifically, the following tables are provided in Appendix~\ref{app:tables}:
\begin{itemize}
\item
Size of the tests under the null hypothesis $H_0$:
\begin{compactitem}
\item Table~\ref{H0main}: Percentage of false rejections when data are serially independent.
\item Table~\ref{H0smmax}: Percentage of false rejections when data are serially dependent.
\end{compactitem}
\item
Power of the tests against specific alternatives:
\begin{compactitem}
\item Table~\ref{H1onecopd2}: Power against $H_{0,m} \cap H_{1,c}$ involving a change of the copula parameter within a copula family and at serial independence.
\item Table~\ref{H1both}: Power against $H_{0,m} \cap H_{1,c}$ involving a change of copula family at a constant value of Kendall's tau and at serial independence.
\item Table~\ref{H0break}: Power against $H_{1,m} \cap H_{0,c}$ involving a change in one of the margins and at serial independence..
\item Table~\ref{H1GHsmmax}: Power against $\neg H_0$ when data are serially dependent and the change occurs in the copula of the innovations.
\end{compactitem}
\end{itemize}

Besides findings of a more anecdotical nature, the following conclusions may be drawn from the results:
\begin{itemize}
\item
All tests hold their level reasonably well in the case of serial independence (Table~\ref{H0main}), with minor fluctuations depending on sample size, test statistic, copula parameter and copula family. 
\item 
In case of serial dependence, the test based on $\hat{S}_n$ is too conservative for the sample sizes under consideration (Table~\ref{H0smmax}). In line with this observation, the test based on $\check{S}_n$ appears to be more powerful than the one based on $\hat S_n$ (Table~\ref{H1GHsmmax}).
\item
For alternative hypotheses involving a change in the \emph{copula}, the tests based on $\hat S_n$ and $\check S_n$ have a higher power than $S_n^R$ (Tables~\ref{H1onecopd2} and~\ref{H1both}). When the copula changes in such a way that Kendall's tau remains constant, the power of $S_n^R$ is especially low. With respect to that last setting, note that distinguishing copulas on the basis of low amounts of data is known to be difficult \citep{GenRemBea09,KojYanHol11}. The fact that the change-point is unknown makes the problem even harder.
\item
For alternative hypotheses involving a change in one of the \emph{margins}, it is the test statistic $S_n^R$ that is substantially more powerful than $S_n$ (Table~\ref{H0break}). The weak power of $S_n$ can be explained by the fact that it is designed for detecting changes in the copula. Another tentative reading of the results is that the test based on $S_n$, regarded as a procedure for testing $H_{0,c}$, is relatively robust against small changes in one margin. In contrast, the test based on $S_n^R$ behaves as an all-purpose test for the hypothesis of a constant distribution rather than as a test for a constant copula.
\end{itemize}

\section{Case studies}
% ====================
\label{sec:cases}

As an illustration, we first applied the test based on $\check S_n$ to bivariate financial data consisting of daily logreturns computed from the DAX and the Standard and Poor 500 indices. Following \citet[Section~7]{DehVogWenWie13}, attention was restricted to the years 2006--2009. The corresponding closing quotes were obtained from \url{http://quote.yahoo.com} using the {\tt get.hist.quote} function of the \texttt{tseries} \textsf{R} package \citep{tseries}, which resulted in $n=993$ bivariate logreturns. Dependent multiplier sequences were generated as explained in Appendix~\ref{strongmixing}. An approximate p-value of 0.04 was obtained, providing some evidence against $H_0$. The conclusion is in line with the results reported in \cite{DehVogWenWie13}. Of course, as discussed earlier, it is only under the assumption that $H_{0,m}$ in~\eqref{H0m} holds that it would be fully justified to decide to reject $H_{0,c}$ in~\eqref{H1c} on the basis of the previous approximate p-value. The value of the change-
point estimator $k_n^\star$ in~\eqref{eq:kstarn} is 529, corresponding to February 22nd, 2008.

As a second illustration, we followed again \cite{DehVogWenWie13} and considered $n=504$ bivariate logreturns computed from closing daily quotes of the Dow Jones Industrial Average and the Nasdaq Composite for the years 1987 and 1988. The former quotes, not being available on \url{http://quote.yahoo.com} anymore, were taken from the \textsf{R} package {\tt QRM} \citep{QRM}. This two-year period is of interest because it contains October 19th, 1987, known as ``black Monday'' \citep[see][Figure~4]{DehVogWenWie13}. An approximate p-value of 0.59 was obtained. Hence, despite the extreme events that occurred during the period under consideration, the test based on $\check S_n$ detects no evidence against $H_0$ in the data, which is in line with the results reported in \cite{DehVogWenWie13}. 

% \ab{With respect to the latter application: it would be interested to test for a break precisely on Black Monday. The definition of an appropriate test statistic is straight-forward. We also investigated this data set in \cite{BucJasWie13} with respect to constancy of the tail dependence (where we applied the tests to (almost independent) standardized residuals from marginal time series models \`a la Remillard and where the hypothesis of constancy over the whole time period got rejected).} \ik{That's interesting and doable. Do you suggest to have it here? I feel very laaaaaazzzzzyyy ;-) .}

\section{Conclusion}
% ==================
\label{sec:conclusion}

We have demonstrated that the sensitivity of rank-based tests for the null hypothesis of a constant distribution against changes in cross-sectional dependence can be improved if ranks are computed with respect to relevant subsamples. In this way, the test we propose achieves in many cases a higher power than the one proposed in \citet{BucRup13}. The limit distribution of the test statistic under the null hypothesis is unwieldy, but approximate p-values can still be computed via a multiplier resampling scheme. To deal with potential serial dependence, we make use of dependent multiplier sequences, an idea going back to \citet{Buh93} and revisited in \cite{BucKoj13}.

Here are some potential avenues for further research:
\begin{compactitem}[--]
\item
Once the null hypothesis has been rejected, the nature of the nonstationary needs to be investigated further: is there a single change-point or is there more than one? Or maybe the change is gradual rather than sudden? And does the change concern the margins or the copula?
% \item
% In Appendix~\ref{app:changepoint}, the estimation of a single change-point is treated. What if there are multiple change-points?
\item
Can one detect a change in the copula without the hypothesis that the margins are constant?
\item
The procedure is computationally intensive because the ranks have to be recomputed for every $k \in \{1, \ldots, n-1\}$. Efficient algorithms for reutilizing calculations from one value of $k$ to the next one might speed up the computations.  
\end{compactitem}

\section*{Acknowledgments}

The authors are grateful to Mark Holmes, Jean-Fran\c{c}ois Quessy and Martin Ruppert for fruitful discussions, and to an anonymous referee for pointing out that Condition~\ref{cond:noTies} might be dispensed with.

The research by A.\ B\"ucher has been supported in parts by the Collaborative Research Center ``Statistical modeling of nonlinear dynamic processes'' (SFB 823, Project A7) of the German Research Foundation (DFG), which is gratefully acknowledged.

J. Segers gratefully acknowledges funding by contract ``Projet d'Act\-ions de Re\-cher\-che Concert\'ees'' No.\ 12/17-045 of the ``Communaut\'e fran\c{c}aise de Belgique'' and by IAP research network Grant P7/06 of the Belgian government (Belgian Science Policy).

\bibliographystyle{plainnat}
\bibliography{biblio}

\begin{thebibliography}{39}
\providecommand{\natexlab}[1]{#1}
\providecommand{\url}[1]{\texttt{#1}}
\expandafter\ifx\csname urlstyle\endcsname\relax
  \providecommand{\doi}[1]{doi: #1}\else
  \providecommand{\doi}{doi: \begingroup \urlstyle{rm}\Url}\fi

\bibitem[Auestad and {Tj\o stheim}(1990)]{AueTjo90}
B.~Auestad and D.~{Tj\o stheim}.
\newblock Identification of nonlinear time series: {F}irst order
  characterization and order determination.
\newblock \emph{Biometrika}, 77:\penalty0 669--687, 1990.

\bibitem[Bai(1994)]{Bai94}
J.~Bai.
\newblock Weak convergence of the sequential empirical processes of residuals
  in {ARMA} models.
\newblock \emph{The Annals of Statistics}, 22\penalty0 (4):\penalty0
  2051--2061, 1994.

\bibitem[B\"ucher(2013{\natexlab{a}})]{Buc13}
A.~B\"ucher.
\newblock A note on nonparametric estimation of bivariate tail dependence.
\newblock \emph{Statistics and Risk Modeling}, page in press,
  2013{\natexlab{a}}.
\newblock URL
  \url{http://www.ruhr-uni-bochum.de/imperia/md/content/mathematik3/publications/taildep_buecher.pdf}.

\bibitem[B\"ucher(2013{\natexlab{b}})]{Buc13b}
A.~B\"ucher.
\newblock A note on weak convergence of the sequential multivariate empirical
  process under strong mixing.
\newblock \emph{Journal of Theoretical Probability}, page in press,
  2013{\natexlab{b}}.

\bibitem[B\"ucher and Dette(2010)]{BucDet10}
A.~B\"ucher and H.~Dette.
\newblock A note on bootstrap approximations for the empirical copula process.
\newblock \emph{Statistics and Probability Letters}, 80\penalty0
  (23--24):\penalty0 1925--1932, 2010.

\bibitem[B\"ucher and Kojadinovic(2013)]{BucKoj13}
A.~B\"ucher and I.~Kojadinovic.
\newblock A dependent multiplier bootstrap for the sequential empirical copula
  process under strong mixing.
\newblock \emph{arXiv:1306.3930}, 2013.

\bibitem[B\"ucher and Ruppert(2013)]{BucRup13}
A.~B\"ucher and M.~Ruppert.
\newblock Consistent testing for a constant copula under strong mixing based on
  the tapered block multiplier technique.
\newblock \emph{Journal of Multivariate Analysis}, 116:\penalty0 208--229,
  2013.

\bibitem[B\"ucher and Segers(2013)]{BucSeg13}
A.~B\"ucher and J.~Segers.
\newblock Extreme value copula estimation based on block maxima of a
  multivariate stationary time series.
\newblock \emph{arXiv:1311.3060}, 2013.

\bibitem[B\"ucher and Volgushev(2013)]{BucVol13}
A.~B\"ucher and S.~Volgushev.
\newblock Empirical and sequential empirical copula processes under serial
  dependence.
\newblock \emph{Journal of Multivariate Analysis}, 119:\penalty0 61--70, 2013.

\bibitem[B\"uhlmann(1993)]{Buh93}
P.~B\"uhlmann.
\newblock \emph{The blockwise bootstrap in time series and empirical
  processes}.
\newblock PhD thesis, ETH Z\"urich, 1993.
\newblock Diss. ETH No. 10354.

\bibitem[Cs\"org\H{o} and Horv{\'a}th(1997)]{CsoHor97}
M.~Cs\"org\H{o} and L.~Horv{\'a}th.
\newblock \emph{Limit theorems in change-point analysis}.
\newblock Wiley Series in Probability and Statistics. John Wiley \& Sons,
  Chichester, UK, 1997.

\bibitem[Dehling et~al.(2013)Dehling, Vogel, Wendler, and Wied]{DehVogWenWie13}
H.~Dehling, D.~Vogel, M.~Wendler, and D.~Wied.
\newblock An efficient and robust test for a change-point in correlation.
\newblock \emph{arXiv:1203.4871}, 2013.

\bibitem[Genest and Segers(2010)]{GenSeg10}
C.~Genest and J.~Segers.
\newblock On the covariance of the asymptotic empirical copula process.
\newblock \emph{Journal of Multivariate Analysis}, 101:\penalty0 1837--1845,
  2010.

\bibitem[Genest et~al.(1998)Genest, Ghoudi, and Rivest]{GenGhoRiv98}
C.~Genest, K.~Ghoudi, and L.-P. Rivest.
\newblock Discussion of ``{U}nderstanding relationships using copulas'', by
  {E}. {F}rees and {E}. {V}aldez.
\newblock \emph{North American Actuarial Journal}, 3:\penalty0 143--149, 1998.

\bibitem[Genest et~al.(2009)Genest, R\'emillard, and Beaudoin]{GenRemBea09}
C.~Genest, B.~R\'emillard, and D.~Beaudoin.
\newblock Goodness-of-fit tests for copulas: {A} review and a power study.
\newblock \emph{Insurance: Mathematics and Economics}, 44:\penalty0 199--213,
  2009.

\bibitem[Gombay and Horv\'ath(1999)]{GomHor99}
E.~Gombay and L.~Horv\'ath.
\newblock Change-points and bootstrap.
\newblock \emph{Environmetrics}, 10\penalty0 (6), 1999.

\bibitem[Gombay and Horv\'ath(2002)]{GomHor02}
E.~Gombay and L.~Horv\'ath.
\newblock Rates of convergence for {$U$}-statistic processes and their
  bootstrapped versions.
\newblock \emph{Journal of Statistical Planning and Inference}, 102:\penalty0
  247--272, 2002.

\bibitem[Hofert et~al.(2013)Hofert, Kojadinovic, M\"achler, and Yan]{copula}
M.~Hofert, I.~Kojadinovic, M.~M\"achler, and J.~Yan.
\newblock \emph{copula: {M}ultivariate dependence with copulas}, 2013.
\newblock URL \url{http://CRAN.R-project.org/package=copula}.
\newblock {R} package version 0.999-7.

\bibitem[Holmes et~al.(2013)Holmes, Kojadinovic, and Quessy]{HolKojQue13}
M.~Holmes, I.~Kojadinovic, and J-F. Quessy.
\newblock Nonparametric tests for change-point detection \`a la {G}ombay and
  {H}orv\'ath.
\newblock \emph{Journal of Multivariate Analysis}, 115:\penalty0 16--32, 2013.

\bibitem[Inoue(2001)]{Ino01}
A.~Inoue.
\newblock Testing for distributional change in time series.
\newblock \emph{Econometric Theory}, 17\penalty0 (1):\penalty0 156--187, 2001.

\bibitem[Khoudraji(1995)]{Kho95}
A.~Khoudraji.
\newblock \emph{Contributions \`a l'\'etude des copules et \`a la
  mod\'elisation des valeurs extr\^emes bivari\'ees}.
\newblock PhD thesis, Universit\'e Laval, Qu\'ebec, Canada, 1995.

\bibitem[Kojadinovic(2014)]{npcp}
I.~Kojadinovic.
\newblock \emph{npcp: {S}ome nonparametric tests for change-point detection in
  (multivariate) observations}, 2014.
\newblock URL \url{http://CRAN.R-project.org/package=npcp}.
\newblock R package version 0.0-1.

\bibitem[Kojadinovic et~al.(2011{\natexlab{a}})Kojadinovic, Segers, and
  Yan]{KojSegYan11}
I.~Kojadinovic, J.~Segers, and J.~Yan.
\newblock Large-sample tests of extreme-value dependence for multivariate
  copulas.
\newblock \emph{The Canadian Journal of Statistics}, 39\penalty0 (4):\penalty0
  703--720, 2011{\natexlab{a}}.

\bibitem[Kojadinovic et~al.(2011{\natexlab{b}})Kojadinovic, Yan, and
  Holmes]{KojYanHol11}
I.~Kojadinovic, J.~Yan, and M.~Holmes.
\newblock Fast large-sample goodness-of-fit for copulas.
\newblock \emph{Statistica Sinica}, 21\penalty0 (2):\penalty0 841--871,
  2011{\natexlab{b}}.

\bibitem[Liebscher(2008)]{Lie08}
E.~Liebscher.
\newblock Construction of asymmetric multivariate copulas.
\newblock \emph{Journal of Multivariate Analysis}, 99:\penalty0 2234--2250,
  2008.

\bibitem[Paparoditis and Politis(2001)]{PapPol01}
E.~Paparoditis and D.N. Politis.
\newblock Tapered block bootstrap.
\newblock \emph{Biometrika}, 88\penalty0 (4):\penalty0 1105--1119, 2001.

\bibitem[Pfaff and McNeil(2013)]{QRM}
B.~Pfaff and A.~McNeil.
\newblock \emph{{QRM: Provides R-language Code to Examine Quantitative Risk
  Management Concepts}}, 2013.
\newblock URL \url{http://CRAN.R-project.org/package=QRM}.
\newblock R package version 0.4-9.

\bibitem[Quessy et~al.(2013)Quessy, Sa\"id, and Favre]{QueSaiFav13}
J.-F. Quessy, M.~Sa\"id, and A.-C. Favre.
\newblock Multivariate {K}endall's tau for change-point detection in copulas.
\newblock \emph{The Canadian Journal of Statistics}, 41:\penalty0 65--82, 2013.

\bibitem[{R Development Core Team}(2013)]{Rsystem}
{R Development Core Team}.
\newblock \emph{{R}: {A} Language and Environment for Statistical Computing}.
\newblock R Foundation for Statistical Computing, Vienna, Austria, 2013.
\newblock URL \url{http://www.R-project.org}.
\newblock {ISBN} 3-900051-07-0.

\bibitem[R\'emillard(2010)]{Rem10}
B.~R\'emillard.
\newblock Goodness-of-fit tests for copulas of multivariate time series.
\newblock \emph{Social Science Research Network}, 1729982:\penalty0 1--32,
  2010.
\newblock URL \url{http://ssrn.com/abstract=1729982}.

\bibitem[R\'emillard and Scaillet(2009)]{RemSca09}
B.~R\'emillard and O.~Scaillet.
\newblock Testing for equality between two copulas.
\newblock \emph{Journal of Multivariate Analysis}, 100\penalty0 (3):\penalty0
  377--386, 2009.

\bibitem[Romano and Wolf(2000)]{RomWol00}
J.P. Romano and M.~Wolf.
\newblock A more general central limit theorem for $m$-dependent random
  variables with unbounded $m$.
\newblock \emph{Statistics and Probability Letters}, 47:\penalty0 115--124,
  2000.

\bibitem[Scaillet(2005)]{Sca05}
O.~Scaillet.
\newblock A {K}olmogorov-{S}mirnov type test for positive quadrant dependence.
\newblock \emph{Canadian Journal of Statistics}, 33:\penalty0 415--427, 2005.

\bibitem[Segers(2012)]{Seg12}
J.~Segers.
\newblock Asymptotics of empirical copula processes under nonrestrictive
  smoothness assumptions.
\newblock \emph{Bernoulli}, 18:\penalty0 764--782, 2012.

\bibitem[Sklar(1959)]{Skl59}
A.~Sklar.
\newblock Fonctions de r\'epartition \`a $n$ dimensions et leurs marges.
\newblock \emph{Publications de l'Institut de Statistique de l'Universit\'e de
  Paris}, 8:\penalty0 229--231, 1959.

\bibitem[Trapletti and Hornik(2013)]{tseries}
A.~Trapletti and K.~Hornik.
\newblock \emph{tseries: {T}ime series analysis and computational finance},
  2013.
\newblock URL \url{http://CRAN.R-project.org/package=tseries}.
\newblock R package version 0.10-32.

\bibitem[{van der Vaart} and Wellner(2000)]{vanWel96}
A.W. {van der Vaart} and J.A. Wellner.
\newblock \emph{Weak convergence and empirical processes}.
\newblock Springer, New York, 2000.
\newblock Second edition.

\bibitem[{van Kampen} and Wied(2013)]{vanWie13}
M.~{van Kampen} and D.~Wied.
\newblock A nonparametric constancy test for copulas under mixing conditions.
\newblock \emph{working paper}, 2013.
\newblock URL \url{http://www.statistik.tu-dortmund.de/vankampen-eng.html}.

\bibitem[Wied et~al.(2013)Wied, Dehling, {van Kampen}, and
  Vogel]{WieDehvanVog13}
D.~Wied, H.~Dehling, M.~{van Kampen}, and D.~Vogel.
\newblock A fluctuation test for constant {S}pearman's rho with nuisance-free
  limit distribution.
\newblock \emph{Computational Statistics and Data Analysis}, page in press,
  2013.

\end{thebibliography}

\appendix

\section{Proof of Proposition~\ref{prop:weak_Cn_sm}}
% =============================================
\label{proof:prop:weak_Cn_sm}

Let us first introduce additional notation. For integers $1 \le k \le l \leq n$, let $H_{k:l}$ denote the empirical c.d.f.\ of the unobservable sample $\vec U_{k}, \dots, \vec U_{l}$ and let $H_{k:l,j}$, for  $j \in \{1, \ldots, d\}$, denote its margins. The empirical quantile functions are
\[
  H_{k:l,j}^{-1}(u) = \inf \{ v \in [0, 1] : H_{k:l,j}(v) \ge u \}, \qquad u \in [0, 1],
\]
which are collected in a vector via
\[
  \vec{H}_{k:l}^{-1}(\vec{u}) = 
  \bigl( H_{k:l,1}^{-1}(u_1), \ldots, H_{k:l,d}^{-1}(u_d) \bigr),
  \qquad \vec{u} \in [0, 1]^d.
\]
By convention, the previously defined quantities are all taken equal to zero if $k > l$.

From the proof of Theorem~3.4 in \cite{BucKoj13}, we have that~\eqref{eq:asymequivCn} holds with $\Cb_n$ replaced by $\Cb_n^{\mathrm{alt}}$, where 
$$
\Cb_n^{\mathrm{alt}}(s,t,\vec u) = \frac{1}{\sqrt{n}} \sum_{i=\ip{ns}+1}^{\ip{nt}}  [ \1 \{ \vec{U}_i \le \vec{H}_{\ip{ns}+1:\ip{nt}}^{-1}(\vec{u}) \} - C(\vec u) ], \qquad (s,t,\vec u) \in \Delta \times [0,1]^d. 
$$
To show~\eqref{eq:asymequivCn}, it remains therefore to prove that 
\begin{equation}
\label{eq:Cbn-Cbnalt}
\sup_{(s,t,\vec u) \in \Delta \times [0,1]^d} \left| \Cb_n^{\mathrm{alt}}(s,t,\vec u) - \Cb_n(s,t,\vec u) \right| \p 0.
\end{equation}
To do so, we adapt the arguments used in Lemma~A.2 of \cite{BucSeg13}. Fix $1 \leq k \leq l \leq n$ and $\vec u \in [0,1]^d$. For $i \in \{ k, \ldots, l \}$, the $d$ components of $\hat{\vec{U}}_i^{k:l}$ defined in~\eqref{eq:pseudo} can be expressed as $\hat{U}_{ij}^{k:l} = H_{k:l,j}(U_{ij})$, $j \in \{ 1, \ldots, d \}$. Next, notice that
\begin{align*}
\1\bigl\{ U_{ij} \le H_{k:l,j}^{-1}(u_j) \bigr\} -  \1\bigl( \hat{U}_{ij}^{k:l} \le u_j \bigr) 
&= \1\bigl\{ U_{ij} < H_{k:l,j}^{-1}(u_j) \bigr\} + \1\bigl\{ U_{ij} = H_{k:l,j}^{-1}(u_j) \bigr\} \\ &\qquad \mbox{} -  \1\bigl( \hat{U}_{ij}^{k:l} < u_j \bigr) -  \1\bigl( \hat{U}_{ij}^{k:l} = u_j \bigr) \\
&= \1\bigl\{ U_{ij} = H_{k:l,j}^{-1}(u_j) \bigr\} - \1\bigl( \hat{U}_{ij}^{k:l} = u_j \bigr), 
\end{align*}
as $x < H^{-1}(u)$ if and only $H(x) < u$ for any distribution function $H$. Since $\hat{U}_{ij}^{k:l} = H_{k:l,j}(U_{ij}) = u_j$ implies $U_{ij} = H_{k:l,j}^{-1}(u_j)$, we obtain that
$$
0 \leq \1\bigl\{ U_{ij} \le H_{k:l,j}^{-1}(u_j) \bigr\} -  \1\bigl( \hat{U}_{ij}^{k:l} \le u_j \bigr) \leq \1\bigl\{ U_{ij} = H_{k:l,j}^{-1}(u_j) \bigr\}.
$$
Combining the previous inequality with the decomposition
\begin{multline*}
\1 \{ \vec U_i \leq \vec{H}_{k:l}^{-1}(\vec{u}) \} - \1 ( \hat{\vec{U}}_i^{k:l} \le \vec{u} ) \\ = \sum_{p=1}^d \left[  \prod_{1 \leq j \leq p}  \1 \{ U_{ij} \leq H_{k:l,j}^{-1}(u_j) \}\prod_{p < j \leq d} \1 ( \hat{U}_{ij}^{k:l} \le u_j ) \right.  \\  \left. - \prod_{1 \leq j \leq p-1}  \1 \{ U_{ij} \leq H_{k:l,j}^{-1}(u_j) \}\prod_{p-1 < j \leq d} \1 ( \hat{U}_{ij}^{k:l} \le u_j ) \right],
\end{multline*}
we obtain that
\begin{equation}
\label{eq:indsineq}
0 \leq \1 \{ \vec U_i \leq \vec{H}_{k:l}^{-1}(\vec{u}) \} - \1 ( \hat{\vec{U}}_i^{k:l} \le \vec{u} ) \leq \sum_{j=1}^d \1\bigl\{ U_{ij} = H_{k:l,j}^{-1}(u_j) \bigr\} . 
\end{equation}
It follows that the supremum in \eqref{eq:Cbn-Cbnalt} is smaller than
$$
\sum_{j=1}^d \sup_{(s,t,u) \in \Delta \times [0,1]} \frac{1}{\sqrt{n}} \sum_{i=\ip{ns}+1}^{\ip{nt}}  \1\bigl\{ U_{ij} = H_{\ip{ns}+1:\ip{nt},j}^{-1}(u) \bigr\} \leq \sum_{j=1}^d \sup_{u \in [0,1]} \frac{1}{\sqrt{n}} \sum_{i=1}^n  \1( U_{ij} = u).
$$
Using the fact that $\1\bigl( U_{ij} = u \bigr) \leq \1\bigl( U_{ij} \leq u \bigr) - \1\bigl( U_{ij} \leq u - 1/n \bigr)$, the latter is smaller 
$$
d \sup_{\substack{\vec u,\vec v \in [0,1]^d \\ \|\vec u - \vec v\|_1 \leq n^{-1}}} \left|  \Bseq(0,1,\vec u) - \Bseq(0,1,\vec v)\right| + d n^{-1/2},
$$
where $\Bseq$ is defined in~\eqref{eq:seqep}. Using the asymptotic uniform equicontinuity in probability of $\B_n$ established in Lemma~2 of \cite{Buc13b}, we finally obtain~\eqref{eq:Cbn-Cbnalt}, which completes the proof. \qed

\section{Proof of Proposition~\ref{prop:check}}
% =============================================
\label{proof:prop:check}

We shall only prove the result in the case of strongly mixing observations, that is, when Condition~\ref{cond:DGP}(ii) is assumed. The proof is similar but simpler when Condition~\ref{cond:DGP}(i) is assumed instead.

It is sufficient to show the statement involving $\check{\Cb}_n^{(m)}$. The statements for $\check{\D}_n^{(m)}$ and $\check{S}_n^{(m)}$ then follow from the continuous mapping theorem.

For any $m \in \{1, \dots, M\}$ and $(s, t, \vec{u}) \in \Delta \times [0, 1]^d$, put
\begin{align}
\label{eq:tildeBnm}
  \Bseq^{(m)} (s, t, \vec{u}) 
  &= \frac{1}{\sqrt{n}} \sum_{i=\ip{ns}+1}^{\ip{nt}} \xi_{i,n}^{(m)} \{\1(\vec U_i \leq \vec{u}) - C(\vec{u}) \}, \\
  \Cb_n^{(m)} (s, t, \vec{u})
\nonumber
  &= \Bseq^{(m)}(s, t, \vec{u}) - \sum_{j=1}^d \dot{C}_j(\vec{u}) \, \Bseq^{(m)}(s, t, \vec{u}^{(j)}).
\end{align}
[Recall that $\vec{u}^{(j)} = (1, \ldots, 1, u_j, 1, \ldots, 1) \in [0,1]^d$, with $u_j$ appearing at the $j$-th coordinate.] From Theorem~2.1 in \cite{BucKoj13}, we have that
\[
  \bigl( \Bseq, \Bseq^{(1)}, \ldots, \Bseq^{(M)} \bigr)
  \leadsto
  \bigl( \B_C, \B_C^{(1)}, \ldots, \B_C^{(M)} \bigr)
\]
in $\{ \ell^\infty( \Delta \times [0, 1]^d ) \}^{M+1}$, where $\B_C^{(1)}, \ldots, \B_C^{(M)}$ are independent copies of $\B_C$ in~\eqref{eq:BC}, and thus, from the continuous mapping theorem and~\eqref{eq:asymequivCn}, we find that
\[
  \bigl( \Cb_n, \Cb_n^{(1)}, \ldots, \Cb_n^{(M)} \bigr)
  \leadsto
  \bigl( \Cb_C, \Cb_C^{(1)}, \ldots, \Cb_C^{(M)} \bigr)
\]
in $\{ \ell^\infty( \Delta \times [0, 1]^d ) \}^{M+1}$. It is therefore sufficient to show that
\begin{equation}
\label{eq:checkCb:conv}
  \sup_{(s, t, \vec{u}) \in \Delta \times [0, 1]^d}
  \bigl| (\check{\Cb}_n^{(m)} - \Cb_n^{(m)}) (s, t, \vec{u}) \bigr|
  \p 0
\end{equation}
for every $m \in \{ 1, \ldots, M \}$. Below, we will show the following two assertions: first,
\begin{equation}
\label{eq:checkBnm:conv}
  \sup_{(s, t, \vec{u}) \in \Delta \times [0, 1]^d}
  \bigl| (\check{\B}_n^{(m)} - \Bseq^{(m)}) (s, t, \vec{u}) \bigr|
  \p 0,
\end{equation}
and second, for every $\delta \in (0, 1/2)$ and every $\eps \in (0, 1)$,
\begin{equation}
\label{eq:dotCjkl:conv}
  \sup_{\substack{\vec{u} \in [0, 1]^d\\\delta \le u_j \le 1-\delta}}
  \sup_{\substack{(s,t) \in [0, 1]^2\\t - s \ge \eps}}
  \bigl|
    \dot{C}_{j,\ip{ns}+1:\ip{nt}}(\vec{u}) - \dot{C}_j(\vec{u})
  \bigr|
  \p 0.
\end{equation}
In view of the structure of $\check{\Cb}_n^{(m)}$ in~\eqref{eq:checkCbnm}, the assertions~\eqref{eq:checkBnm:conv} and~\eqref{eq:dotCjkl:conv} imply~\eqref{eq:checkCb:conv}, as we show next. Clearly,
\begin{multline}
\label{eq:threeterms}
  \bigl| (\check{\Cb}_n^{(m)} - \Cb_n^{(m)}) (s, t, \vec{u}) \bigr| \\
  \le \bigl| (\check{\B}_n^{(m)} - \Bseq^{(m)}) (s, t, \vec{u}) \bigr|
  + \sum_{j=1}^d \bigl| \dot C_{j,\ip{ns}+1:\ip{nt}}(\vec{u}) \bigr| \,
    \bigl| (\check{\B}_n^{(m)} - \Bseq^{(m)}) (s, t, \vec{u}^{(j)}) \bigr| \\
  + \sum_{j=1}^d 
  \bigl| \dot C_{j,\ip{ns}+1:\ip{nt}}(\vec{u}) - \dot{C}_j(\vec{u}) \bigr| \, 
  \bigl| \Bseq^{(m)}(s, t, \vec{u}^{(j)}) \bigr|.
\end{multline}
Taking suprema over $(s, t, \vec{u}) \in \Delta \times [0, 1]^d$, the first and the second term on the right-hand side of~\eqref{eq:threeterms} converge to zero in probability because of assertion~\eqref{eq:checkBnm:conv} and uniform boundedness of $\dot{C}_{j,k:l}$ \citep[see][proof of Proposition~2]{KojSegYan11}. The third term on the right-hand side of~\eqref{eq:threeterms} converges to zero in probability because of assertion~\eqref{eq:dotCjkl:conv} and the fact that $(s, t, \vec{u}) \mapsto \Bseq^{(m)}(s, t, \vec{u}^{(j)})$ vanishes as soon as $s = t$ or $u_j \in \{ 0, 1 \}$, and is asymptotically uniformly equicontinuous in probability as a consequence of Lemma~A.3 in \cite{BucKoj13}.

It remains to show~\eqref{eq:checkBnm:conv} and~\eqref{eq:dotCjkl:conv}. The proof of the latter assertion is simplest and is given first.

\begin{proof}[Proof of~\eqref{eq:dotCjkl:conv}]
Observe that
\[
  C_{\ip{ns}+1:\ip{nt}}( \vec{u} )
  = C( \vec{u} ) + \frac{1}{\sqrt{n} \, \lambda_n(s, t)} \, \Cb_n(s, t, \vec{u}).
\]
Fix $\delta \in (0,1/2)$ and $\eps \in (0,1)$. Without loss of generality, assume that $n$ is large enough so that the bandwidth $h = h_n(s, t) = 1/\sqrt{\ip{nt} - \ip{ns}}$ is less than $\delta$ whenever $t - s \ge \eps$. Then, for $t - s \ge \eps$ and $\vec u \in [0, 1]^d$ with $\delta \le u_j \le 1-\delta$, we have
\begin{multline*}
  \dot{C}_{j,\ip{ns}+1:\ip{nt}}(\vec{u})
  = \frac{1}{2h} \{ C( \vec{u} + h \vec{e}_j ) - C( \vec{u} - h \vec{e}_j ) \} \\
  + \frac{1}{2h \, \sqrt{n} \, \lambda_n(s, t)} \{ \Cb_n(s, t, \vec{u} + h \vec{e}_j) - \Cb_n(s, t, \vec{u} - h \vec{e}_j) \}.
\end{multline*}
By the assumption of existence and continuity of $\dot{C}_j$ on $V_j$ (see Condition~\ref{cond:pd}), and since $0 \le \dot{C}_j \le 1$, it follows from the mean-value theorem that
\[
  \sup_{\substack{\vec{u} \in [0, 1]^d\\\delta \le u_j \le 1-\delta}}
  \biggl|
    \frac{1}{2h} \{ C( \vec{u} + h \vec{e}_j ) - C( \vec{u} - h \vec{e}_j ) \}
    - \dot{C}_j(\vec{u})
  \biggr|
  \to 0, \qquad h \to 0.
\]
Using~\eqref{eq:asymequivCn} %uniform boundedness of $\dot C_j$ (just above ;-) )
and the fact that $\Bseq$ is asymptotically uniformly equicontinuous in probability, %as a consequence of Lemma~2 in \cite{Buc13b}, 
it can be verified that $\Cb_n$ is asymptotically uniformly equicontinuous in probability as well. It follows that
\[
  \sup_{\substack{\vec{u} \in [0, 1]^d\\\delta \le u_j \le 1-\delta}}
  \sup_{\substack{(s,t) \in [0, 1]^2\\t - s \ge \eps}}
  \bigl| 
    \Cb_n(s, t, \vec{u} + h \vec{e}_j) - \Cb_n(s, t, \vec{u} - h \vec{e}_j) 
  \bigr|
  \p 0.
\]
Finally,
\[
  \frac{1}{2h \, \sqrt{n} \, \lambda_n(s, t)}
  = \frac{1}{2 \sqrt{ \lambda_n(s, t) }} \le \frac{1}{2 \sqrt{\eps - 1/n}}.
\]
Combine the four previous displays to arrive at the desired conclusion.
\end{proof}

The proof of~\eqref{eq:checkBnm:conv} is more complicated. Using the notation introduced in Appendix~\ref{proof:prop:weak_Cn_sm}, let us define the auto-centered version of the process $\Bseq^{(m)}$ in~\eqref{eq:tildeBnm} as 
\begin{align}
\nonumber
  \Bseqcen^{(m)} (s, t, \vec{u}) 
  &= \frac{1}{\sqrt{n}} \sum_{i=\ip{ns}+1}^{\ip{nt}} \xi_{i,n}^{(m)} \{\1(\vec U_i \leq \vec{u}) - H_{\ip{ns}+1:\ip{nt}}(\vec{u}) \} \\
  \label{eq:checktildeBnm}
  &= \frac{1}{\sqrt{n}} \sum_{i=\ip{ns}+1}^{\ip{nt}} ( \xi_{i,n}^{(m)} - \bar  \xi_{\ip{ns}+1:\ip{nt}}^{(m)} ) \, \1 ( \vec U_i \leq \vec{u} ),
\end{align}
with the usual convention that empty sums are zero.

\begin{proof}[Proof of~\eqref{eq:checkBnm:conv}]
Consider the decomposition
\begin{align*}
  \bigl| (\check{\B}_n^{(m)} - \Bseq^{(m)}) (s, t, \vec{u}) \bigr| 
  &\le
  \bigl| 
    \check{\B}_n^{(m)} (s, t, \vec{u}) 
    - 
    \Bseqcen^{(m)} (s, t, \vec{H}_{\ip{ns}+1:\ip{nt}}^{-1}(\vec{u})) 
  \bigr| \\
  &\phantom{\le} + 
  \bigl| 
    \Bseqcen^{(m)} (s, t, \vec{H}_{\ip{ns}+1:\ip{nt}}^{-1}(\vec{u})) 
    - 
    \Bseq^{(m)} (s, t, \vec{H}_{\ip{ns}+1:\ip{nt}}^{-1}(\vec{u})) 
  \bigr| \\
  &\phantom{\le} +
  \bigl|
    \Bseq^{(m)} (s, t, \vec{H}_{\ip{ns}+1:\ip{nt}}^{-1}(\vec{u})) 
    -
    \Bseq^{(m)} (s, t, \vec{u})
  \bigr|.
\end{align*}
Write out the definitions of the processes $\Bseq^{(m)}$, $\Bseqcen^{(m)}$ and $\check{\B}_n^{(m)}$ in~\eqref{eq:tildeBnm},~\eqref{eq:checktildeBnm} and~\eqref{eq:checkBnm}, respectively, and take suprema over $(s, t, \vec{u}) \in \Delta \times [0, 1]^d$ to obtain
\begin{align}
    \sup_{(s, t, \vec{u}) \in \Delta \times [0, 1]^d} 
    &\bigl| 
      (\check{\B}_n^{(m)} - \Bseq^{(m)}) (s, t, \vec{u}) 
    \bigr|
   \nonumber \\
  &\le
  \sup_{(s, t, \vec{u}) \in \Delta \times [0, 1]^d}  \bigl| 
    \check{\B}_n^{(m)} (s, t, \vec{u}) 
    - 
    \Bseqcen^{(m)} (s, t, \vec{H}_{\ip{ns}+1:\ip{nt}}^{-1}(\vec{u})) 
  \bigr|
\label{eq:termB:1} \\
  &\phantom{\le} +
  \sup_{(s, t, \vec{u}) \in \Delta \times [0, 1]^d}
  \biggl|
    \frac{1}{\sqrt{n}} \sum_{i=\ip{ns}+1}^{\ip{nt}} \xi_{i,n}^{(m)} \,     
  \biggr| \,
  \bigl| H_{\ip{ns}+1:\ip{nt}}(\vec{u}) - C(\vec{u}) \bigr|
  \label{eq:termB:2} \\
  &\phantom{\le} +
  \sup_{(s, t, \vec{u}) \in \Delta \times [0, 1]^d}
  \bigl|
    \Bseq^{(m)} (s, t, \vec{H}_{\ip{ns}+1:\ip{nt}}^{-1}(\vec{u})) 
    -
    \Bseq^{(m)} (s, t, \vec{u})
  \bigr|. \label{eq:termB:3}
\end{align}
Each term requires a different treatment.
\begin{enumerate}[1-]
% \item \emph{The term~\eqref{eq:termB:1}:}
% Fix $1 \leq k \leq l \leq n$. Writing $\vec{U}_i = (U_{i1}, \ldots, U_{id})$, the $d$ components of $\hat{\vec{U}}_i^{k:l}$ can be expressed as
% \[
%   \hat{U}_{ij}^{k:l} = H_{k:l,j}(U_{ij}), \qquad i \in \{ k, \ldots, l \}.
% \]
% The indicators $\1( \hat{U}_{ij}^{k:l} \le u_j )$ and $\1( U_{ij} \le H_{k:l,j}^{-1}(u_j) )$ can be different only if $U_{ij} = H_{k:l,j}^{-1}(u_j)$. Because of the ``no-ties'' Condition~\ref{cond:noTies} and the continuity of the marginal distributions, for every $j \in \{1, \ldots, d\}$, there is at most one index $i \in \{ k, \ldots, l \}$ for which the indicators are different. It follows that the term~\eqref{eq:termB:1} is bounded by
% \[
%   \frac{2d}{\sqrt{n}} \max_{1 \leq i \leq n} | \xi_{i,n}^{(m)} |.
% \]
% Using the fact that, from~(M1), for any $\nu \geq 1$, $\sup_{n \geq 1} \expec[ | \xi_{1,n}^{(m)} |^\nu ] < \infty$, we have that, for every $\alpha > 0$ and $\nu \geq 1$ such that $\nu > 1 / \alpha$,
% \begin{equation}
% \label{eq:maxxi}
%   \Pr \biggl( \max_{1 \leq i \leq n} | \xi_{i,n}^{(m)} | \ge n^{\alpha} \biggr)
%   \le n \Pr \bigl( | \xi_{1,n}^{(m)} | \ge n^{\alpha} \bigr)
%   \le n^{1 - \nu \alpha} \, \expec[ | \xi_{1,n}^{(m)} |^\nu ]
%   \to 0.
% \end{equation}
% Take $\alpha \in (0, 1/2)$ to conclude that~\eqref{eq:termB:1} converges to zero in probability.
\item \emph{The term~\eqref{eq:termB:2}:}
Let $\delta \in (0, 1/2)$, to be specified later. We split the supremum into two parts, according to whether $t-s$ is smaller or larger than $a_n = n^{-1/2-\delta}$:
\begin{align*}
  A_{n,1}
  &= \sup_{\substack{(s, t, \vec{u}) \in \Delta \times [0, 1]^d\\t-s \le a_n}}
  \biggl|
    \frac{1}{\sqrt{n}}
    \sum_{i=\ip{ns}+1}^{\ip{nt}} \xi_{i,n}^{(m)}
  \biggr| \;
  \bigl| 
    H_{\ip{ns}+1:\ip{nt}}( \vec{u} ) - C( \vec{u} ) 
  \bigr|, \\
  A_{n,2}
  &= \sup_{\substack{(s, t, \vec{u}) \in \Delta \times [0, 1]^d\\t-s \ge a_n}}
  \biggl|
    \frac{1}{\sqrt{n}}
    \sum_{i=\ip{ns}+1}^{\ip{nt}} \xi_{i,n}^{(m)}
  \biggr| \;
  \bigl| 
    H_{\ip{ns}+1:\ip{nt}}( \vec{u} ) - C( \vec{u} ) 
  \bigr|.
\end{align*}
We will show that both $A_{n,1}$ and $A_{n,2}$ converge to zero in probability.
\begin{enumerate}[(a)]
\item
Since both $H_{k:l}$ and $C$ take values in $[0, 1]$, a crude bound for $A_{n,1}$ is
\[
  A_{n,1} 
  \le \frac{1}{\sqrt{n}} (n \, a_n + 1) \max_{1 \le i \le n} | \xi_{i,n}^{(m)} |
  \le 2 \, n^{-\delta} \, \max_{1 \le i \le n} | \xi_{i,n} |.
\]
Using the fact that, from~(M1), for any $\nu \geq 1$, $\sup_{n \geq 1} \expec[ | \xi_{1,n}^{(m)} |^\nu ] < \infty$, we have that, for every $\alpha > 0$ and $\nu \geq 1$ such that $\nu > 1 / \alpha$,
%\begin{equation}
%\label{eq:maxxi}
$$
\Pr \biggl( \max_{1 \leq i \leq n} | \xi_{i,n}^{(m)} | \ge n^{\alpha} \biggr)
\le n \Pr \bigl( | \xi_{1,n}^{(m)} | \ge n^{\alpha} \bigr)
\le n^{1 - \nu \alpha} \, \expec[ | \xi_{1,n}^{(m)} |^\nu ]
\to 0.
$$
%\end{equation}
Apply the previous display with $\alpha \in (0, \delta)$ to find that $A_{n,1}$ converges to zero in probability.
\item
Recall $\B_n$ in~\eqref{eq:seqep}. Observe that
\[
  \B_n(s, t, \vec{u}) 
  = \frac{\ip{nt} - \ip{ns}}{\sqrt{n}} \, \{ H_{\ip{ns}+1:\ip{nt}}(\vec{u}) - C(\vec{u}) \}.
\]
We have
\begin{align*}
  A_{n,2}
  &= \sup_{\substack{(s, t, \vec{u}) \in \Delta \times [0, 1]^d\\t-s \ge a_n}}
  \biggl|
    \frac{1}{\ip{nt}-\ip{ns}}
    \sum_{i=\ip{ns}+1}^{\ip{nt}} \xi_{i,n}^{(m)}
  \biggr| \;
  \bigl| 
    \B_n(s, t, \vec{u})
  \bigr| \\
 % &\le \sup_{\substack{\ip{ns} < \ip{nt} \\ nt-ns \ge na_n}}
 %  \biggl|
 %    \frac{1}{\ip{nt}-\ip{ns}}
 %    \sum_{i=\ip{ns}+1}^{\ip{nt}} \xi_{i,n}^{(m)}
 %  \biggr| \;
 %  \sup_{(s, t, \vec{u}) \in \Delta \times [0, 1]^d}
 %  \bigl| 
 %    \B_n(s, t, \vec{u})
 %  \bigr| \\ 
  &\le \sup_{\substack{\ip{ns} < \ip{nt}\\ \ip{nt}+1-\ip{ns} \ge na_n}}
  \biggl|
    \frac{1}{\ip{nt}-\ip{ns}}
    \sum_{i=\ip{ns}+1}^{\ip{nt}} \xi_{i,n}^{(m)}
  \biggr| \;
  \sup_{(s, t, \vec{u}) \in \Delta \times [0, 1]^d}
  \bigl| 
    \B_n(s, t, \vec{u})
  \bigr| \\ 
  % &\le \max_{\substack{1 \le k < l \le n\\l-k \ge na_n - 1}}
  % \biggl|
  %   \frac{1}{l-k}
  %   \sum_{i=k+1}^{l} \xi_{i,n}^{(m)}
  % \biggr| \;
  % \sup_{(s, t, \vec{u}) \in \Delta \times [0, 1]^d}
  % \bigl| 
  %   \B_n(s, t, \vec{u})
  % \bigr| \\ 
   &\le \max_{\substack{1 \le k \le l \le n\\l-k \ge na_n}}
  \biggl|
    \frac{1}{l-k+1}
    \sum_{i=k}^{l} \xi_{i,n}^{(m)}
  \biggr| \;
  \sup_{(s, t, \vec{u}) \in \Delta \times [0, 1]^d}
  \bigl| 
    \B_n(s, t, \vec{u})
  \bigr|. 
  % &\le
  % \max_{\substack{1 \le k \le l \le n\\l-k+1 \ge na_n - 1}}
  % \biggl|
  %   \frac{1}{l-k+1}
  %   \sum_{i=k}^l \xi_{i,n}^{(m)}
  % \biggr|
  % \;
  % \sup_{(s, t, \vec{u}) \in \Delta \times [0, 1]^d}
  % \bigl| 
  %   \B_n(s, t, \vec{u})
  % \bigr| \mbox{\ik{previous final inequality}}
\end{align*}
By weak convergence $\Bseq \leadsto \B_C$ in $\ell^\infty(\Delta \times [0, 1]^d)$, the supremum at the end of the previous display is bounded in probability. Writing $b_n = n a_n$, it is sufficient to show that
\[
  \max_{\substack{1 \le k \le l \le n\\ l-k \ge b_n}}
  \biggl|
    \frac{1}{l-k+1} \sum_{i=k}^l \xi_{i,n}^{(m)}
  \biggr|
  \p 0, \qquad n \to \infty.
\]
Fix $\eta > 0$. The probability that the previous maximum exceeds $\eta$ is bounded by
\begin{equation}
\label{eq:prob:sum}
  \sum_{\substack{1 \le k \le l \le n\\\l-k \ge b_n}}
  \Pr 
  \biggl[
    \biggl|
      \frac{1}{l-k+1} \sum_{i=k}^l \xi_{i,n}^{(m)}
    \biggr|
    > \eta
  \biggr].
\end{equation}
Fix $\nu \ge 2$, to be specified later. By stationarity and Markov's inequality, the previous expression is bounded by
\[
  \sum_{\substack{1 \le k \le l \le n\\l-k \ge b_n}}
  \eta^{-\nu} \, (l-k+1)^{-\nu} \, \expec \biggl[ \biggl| \sum_{i=1}^{l-k+1} \xi_{i,n}^{(m)} \biggr|^\nu \biggr]
  \le
  \eta^{-\nu} \, n \sum_{b_n \leq r \leq n} r^{-\nu} \, \expec \biggl[ \biggl| \sum_{i=1}^r \xi_{i,n}^{(m)} \biggr|^\nu \biggr].
\]
Recall that the sequence $(\xi_{i,n}^{(m)})_{i \in \Z}$ is $\ell_n$-dependent from~(M2) and assume that $n$ is sufficiently large so that $n \geq 2 \ell_n + b_n$. Then, by Corollary~A.1 in \cite{RomWol00}, there exists  a constant $C_\nu$, depending only on $\nu$, such that
\[
  \expec \biggl[ \biggl| \sum_{i=1}^r \xi_{i,n}^{(m)} \biggr|^\nu \biggr]
  \le C_\nu^\nu \, (4 \ell_n r)^{\nu/2} \, \expec[ | \xi_{1,n}^{(m)} |^\nu ].
\]
Using the fact that, from~(M1), $\sup_{n\geq1} \expec[ | \xi_{1,n}^{(m)} |^\nu ] < \infty$, and up to a multiplicative constant, the expression in~\eqref{eq:prob:sum} is bounded by
\begin{align*}
  n \sum_{b_n \leq r \leq n} r^{-\nu} (\ell_n r)^{\nu / 2}
  &\le n^2 b_n^{-\nu/2} \ell_n^{\nu / 2} \\
  &= O \bigl( n^{2 - (1/2-\delta)\nu/2 + (1/2 - \gamma) \nu / 2} \bigr) \\
  &= O \bigl( n^{2 + (\delta - \gamma) \nu / 2} \bigr).
\end{align*}
The right-hand side converges to zero if we choose $\delta = \gamma / 2$ and then $\nu > 8 / \gamma$.
\end{enumerate}

\item \emph{The term~\eqref{eq:termB:3}:} We have to show that, for every $\eta,\lambda > 0$,
$$
\Pr \left[ \sup_{(s, t, \vec{u}) \in \Delta \times [0, 1]^d}
  \bigl|
    \Bseq^{(m)} \{ s, t, \vec{H}_{\ip{ns}+1:\ip{nt}}^{-1}(\vec{u}) \}
    -
    \Bseq^{(m)} (s, t, \vec{u}) \bigr| > \lambda \right] \leq \eta,
$$
for all sufficiently large $n$. 

Fix $\eta,\lambda > 0$. Since~\eqref{eq:termB:3} is smaller than $2 \sup_{(s, t, \vec{u}) \in \Delta \times [0, 1]^d} | \Bseq^{(m)} (s, t, \vec{u})|$, and using the fact that $\Bseq^{(m)}$ vanishes on the diagonal $s=t$ and is asymptotically uniformly equicontinuous in probability, %as a consequence of Lemma~A.3 in \cite{BucKoj13}, 
there exists $\eps \in (0,1)$ such that, for all $n$ sufficiently large, 
$$
\Pr \left[ \sup_{\substack{(s, t, \vec{u}) \in \Delta \times [0, 1]^d \\ t-s < \eps}}
  \bigl|
    \Bseq^{(m)} \{ s, t, \vec{H}_{\ip{ns}+1:\ip{nt}}^{-1}(\vec{u}) \}
    -
    \Bseq^{(m)} (s, t, \vec{u}) \bigr| > \lambda \right] \leq \eta/2.
$$
Setting $\zeta_n = \sup_{\substack{(s,t,\vec u) \in \Delta \times [0,1]^d \\ t-s \geq \eps}} \| \vec{H}_{\ip{ns}+1:\ip{nt}}^{-1}(\vec u) - \vec u \|_1$, we shall now show that, for all $n$ sufficiently large, 
$$
A_n = \Pr \left\{ \sup_{\substack{(s, t, \vec u, \vec v) \in \Delta \times [0, 1]^{2d} \\ t-s \geq \eps, \, \|\vec u - \vec v\|_1 \leq \zeta_n}}
  \bigl|
    \Bseq^{(m)} ( s, t, \vec u )  -
    \Bseq^{(m)} (s, t, \vec v) \bigr| > \lambda \right\} \leq \eta/2,
$$
which will complete the proof. Using again the asymptotic uniformly equicontinuity in probability of $\Bseq^{(m)}$, there exists $\mu \in (0,1)$ such that, for all $n$ sufficiently large, 
 $$
A_{n,1} = \Pr \left\{ \sup_{\substack{(s, t, \vec u, \vec v) \in \Delta \times [0, 1]^{2d} \\ t-s \geq \eps, \, \|\vec u - \vec v\|_1 \leq \mu}}
  \bigl|
    \Bseq^{(m)} ( s, t, \vec u )  -
    \Bseq^{(m)} (s, t, \vec v) \bigr| > \lambda \right\} \leq \eta/4.
$$
We then bound $A_n$ by $A_{n,1} + A_{n,2}$, where $A_{n,2} = \Pr(\zeta_n > \mu)$. From the weak convergence of $\Bseq$ to $\B_C$ in $\ell^\infty(\Delta \times [0,1]^d)$, we have that 
\begin{multline*}
\sup_{\substack{(s,t,\vec u) \in \Delta \times [0,1]^d \\ t-s \geq \eps}} |H_{\ip{ns}+1:\ip{nt}}(\vec u) - C(\vec u)| \\ \leq \sup_{(s,t,\vec u) \in \Delta \times [0,1]^d} |\B_n(s,t,\vec u)| \times n^{-1/2} \times \sup_{\substack{(s,t) \in \Delta \\ t-s \geq \eps}} \{\lambda_n(s,t)\}^{-1} \p 0.
\end{multline*}
Using the fact that  $\sup_{u \in [0,1]} |H_{\ip{ns}+1:\ip{nt},j}(u) - u | = \sup_{u \in [0,1]} |H_{\ip{ns}+1:\ip{nt},j}^{-1}(u) - u |$ for $j \in \{1,\dots,d\}$ (for instance, by symmetry arguments on the graphs of $H_{\ip{ns}+1:\ip{nt},j}$ and $H_{\ip{ns}+1:\ip{nt},j}^{-1}$), we immediately obtain that $\zeta_n \p 0$, which implies that, for all $n$ sufficiently large, $A_{n,2} \leq \eta/4$, and thus that, for all $n$ sufficiently large, $A_n \leq \eta/2$.

\item \emph{The term~\eqref{eq:termB:1}:} For the following arguments, it is sufficient to assume that the sequence $(\xi_{i,n}^{(m)})_{i \in \Z}$ appearing in $\check \B_n^{(m)}$ and $\Bseqcen^{(m)}$ satisfies only (M1) with $\Ex[ \{ \xi_{0,n}^{(m)} \}^2 ] > 0$ not necessarily equal to one. 

Let $K > 0$ be a constant and let us first suppose that, for any  $n \geq 1$ and $i \in \{1,\dots,n\}$, $\xi_{i,n}^{(m)} \geq -K$. With~\eqref{eq:indsineq} in mind, the term~\eqref{eq:termB:1} is smaller than $A_{n,1} + A_{n,2}$, where
$$
A_{n,1} = \sup_{(s,t,\vec u) \in \Delta \times [0,1]^d} %\left| 
 \frac{1}{\sqrt{n}} \sum_{i=\ip{ns}+1}^{\ip{nt}} ( \xi_{i,n}^{(m)} + K ) \left[ \1 \{ \vec U_i \leq \vec{H}_{\ip{ns}+1:\ip{nt}}^{-1}(\vec{u}) \} - \1 ( \hat{\vec{U}}_i^{\ip{ns}+1:\ip{nt}} \le \vec{u} )  \right]% \right|
$$
and
$$
A_{n,2} = \sup_{(s,t,\vec u) \in \Delta \times [0,1]^d} %\left|  
\frac{K  + \bar  \xi_{\ip{ns}+1:\ip{nt}}^{(m)}}{\sqrt{n}} \sum_{i=\ip{ns}+1}^{\ip{nt}}   \left[ \1 \{ \vec U_i \leq \vec{H}_{\ip{ns}+1:\ip{nt}}^{-1}(\vec{u}) \} - \1 ( \hat{\vec{U}}_i^{\ip{ns}+1:\ip{nt}} \le \vec{u} )  \right]. % \right|.
$$
Let us first show that $A_{n,1} \p 0$. Plugging~\eqref{eq:indsineq} into the expression of $A_{n,1}$, we bound $A_{n,1}$ by $A_{n,1,1} + \dots + A_{n,1,d}$, where
$$
A_{n,1,j} = \sup_{(s,t,u) \in \Delta \times [0,1]} \frac{1}{\sqrt{n}} \sum_{i=\ip{ns}+1}^{\ip{nt}} ( \xi_{i,n}^{(m)} + K ) \1\bigl( U_{ij} = u \bigr).
$$ 
To prove that $A_{n,1} \p 0$, we shall now show that $A_{n,1,j} \p 0$ for all $j \in \{1,\dots,d\}$. Fix $j \in \{1,\dots,d\}$. Using the fact that $\1( U_{ij} = u ) \leq \1( U_{ij} \leq u ) - \1( U_{ij} \leq u - 1/n)$, we obtain that $A_{n,1,j}$ is smaller than $A_{n,1,j}' + A_{n,1,j}'' + A_{n,1,j}'''$, where
\begin{align*}
A_{n,1,j}' &= \sup_{(s,t,u) \in \Delta \times [0,1]} \left|  \frac{1}{\sqrt{n}} \sum_{i=\ip{ns}+1}^{\ip{nt}} ( \xi_{i,n}^{(m)} - \bar  \xi_{\ip{ns}+1:\ip{nt}}^{(m)}) \left\{ \1\bigl( U_{ij} \leq u \bigr) - \1\bigl( U_{ij} \leq u - 1/n \bigr) \right\}  \right|, \\
A_{n,1,j}'' &= K \sup_{(s,t,u) \in \Delta \times [0,1]} \frac{1}{\sqrt{n}} \sum_{i=\ip{ns}+1}^{\ip{nt}} \left\{ \1\bigl( U_{ij} \leq u \bigr) - \1\bigl( U_{ij} \leq u - 1/n \bigr) \right\} , \\
A_{n,1,j}''' &= \sup_{(s,t,u) \in \Delta \times [0,1]}  \frac{\left| \bar \xi_{\ip{ns}+1:\ip{nt}}^{(m)} \right|}{\sqrt{n}} \sum_{i=\ip{ns}+1}^{\ip{nt}} \left\{ \1\bigl( U_{ij} \leq u \bigr) - \1\bigl( U_{ij} \leq u - 1/n \bigr) \right\}  . 
\end{align*}
From Lemma A.3 of \cite{BucKoj13}, we know that $\Bseq^{(m)}$ is asymptotically uniformly equicontinuous in probability under the weaker conditions on the sequence $(\xi_{i,n}^{(m)})_{i \in \Z}$  considered above. The treatment of the term~\eqref{eq:termB:2} carried out previously remains valid under these conditions and ensures that
$$
\sup_{(s,t,\vec u) \in \Delta \times [0,1]^d}\bigl| 
    \Bseqcen^{(m)} (s, t, \vec{u}) 
    - 
    \Bseq^{(m)} (s, t,\vec{u}) 
  \bigr| \p 0,
$$
which implies that $\Bseqcen^{(m)}$ is asymptotically uniformly equicontinuous in probability as well under the same weaker conditions on the sequence $(\xi_{i,n}^{(m)})_{i \in \Z}$. The latter immediately implies that $A_{n,1,j}' \p 0$. For $A_{n,1,j}''$, we have
$$
A_{n,1,j}'' \leq K \sup_{\substack{\vec u,\vec v \in [0,1]^d \\ \|\vec u - \vec v\|_1 \leq n^{-1}}} \left|  \Bseq(0,1,\vec u) - \Bseq(0,1,\vec v)\right| + K n^{-1/2} \p 0,
$$
by asymptotic uniform equicontinuity in probability of $\Bseq$. The fact that $A_{n,1,j}''' \p 0$ can be shown by proceeding as for the term~\eqref{eq:termB:2}. Hence, we have that $A_{n,1,j} \p 0$, which implies that $A_{n,1} \p 0$. 

The fact that $A_{n,2} \p 0$, follows from~\eqref{eq:indsineq} which implies that $A_{n,2}$ is smaller than $\sum_{j=1}^d (A_{n,1,j}'' + A_{n,1,j}''')$. This completes the proof under the condition $\xi_{i,n}^{(m)} \geq -K$.

To show that this condition is not necessary, we proceed as at the end of the proof of Lemma A.3 of \cite{BucKoj13}. Let $Z_{i,n}^+=\max(\xi_{i,n}^{(m)},0)$, $Z_{i,n}^- = \max(-\xi_{i,n}^{(m)},0)$, $K^+ = \Ex(Z_{0,n}^+)$ and  $K^- = \Ex(Z_{0,n}^-)$. Furthermore, define $\xi_{i,n}^{(m),+} = Z_{i,n}^+- K^+$ and $\xi_{i,n}^{(m),-} = Z_{i,n}^- - K^-$. Then, using the fact that $K^+ - K^- = 0$, we can write
\[
\xi_{i,n}^{(m)}  = Z_{i,n}^+ - Z_{i,n}^- = Z_{i,n}^+- K^+ - ( Z_{i,n}^-- K^- ) = \xi_{i,n}^{(m),+} - \xi_{i,n}^{(m),-}.
\]
Let $\Bseq^{(m),+}$ and $\Bseq^{(m),-}$ be the analogues of $\Bseq^{(m)}$ defined from the sequences $(\xi_{i,n}^{(m),+})_{i \in \Z}$ and $(\xi_{i,n}^{(m),-})_{i \in \Z}$, respectively, and similarly for $\Bseqcen^{(m),+}$ and $\Bseqcen^{(m),-}$. The case treated above yields 
\begin{align*}
&\sup_{(s,t,\vec u) \in \Delta \times [0,1]^d} \bigl| 
    \check{\B}_n^{(m),+} (s, t, \vec{u}) 
    - 
    \Bseqcen^{(m),+} (s, t, \vec{H}_{\ip{ns}+1:\ip{nt}}^{-1}(\vec{u})) 
  \bigr| \p 0, \\
&\sup_{(s,t,\vec u) \in \Delta \times [0,1]^d} \bigl| 
    \check{\B}_n^{(m),-} (s, t, \vec{u}) 
    - 
    \Bseqcen^{(m),-} (s, t, \vec{H}_{\ip{ns}+1:\ip{nt}}^{-1}(\vec{u})) 
  \bigr| \p 0.
\end{align*}
The desired result finally follows from the fact that $\Bseq^{(m)} = \Bseq^{(m),+}  - \Bseq^{(m),-}$ and $\Bseqcen^{(m)} = \Bseqcen^{(m),+}  - \Bseqcen^{(m),-}$.
%%%%%%%%%%%%%%%%%%%%%%%%%%%%%%%%%%%%%%%%%%%%%%%%%%%%%%%%%%%%%%%%%%
\end{enumerate}
\end{proof}

\section{On the set-up of the simulation experiments}
% ========================================================
\label{strongmixing}

For the numerical experiments involving serially dependent observations, we restricted ourselves to the bivariate case and only focused on the tests based on $\check S_n$ and $\hat S_n$. Given a bivariate copula $C$, two models were used to generate serially dependent observations under $H_0$ defined in~\eqref{H0}. 
\begin{itemize}
\item
The first one is a simple autoregressive model of order one, AR(1). Let $\vec U_i$, $i \in \{-100,\dots,0,\dots,n\}$, be a bivariate i.i.d.\ sample from a copula~$C$. Then, set $\vec \epsilon_i = (\Phi^{-1}(U_{i1}),\Phi^{-1}(U_{i2}))$, where $\Phi$ is the c.d.f.\ of the standard normal distribution, and $\vec X_{-100} = \vec \epsilon_{-100}$. Finally, for any $j \in \{1,2\}$ and $i \in \{-99,\dots,0,\dots,n\}$, compute recursively
\begin{equation}
\tag{AR1}
\label{eq:ar1}
  X_{ij} = 0.5 X_{i-1,j} + \epsilon_{ij}.
\end{equation}
\item
The second model is a bivariate version of the exponential autoregressive (EXPAR) model considered in \cite{AueTjo90} and \citet[Section~3.3]{PapPol01} \citep[see also][]{BucKoj13}. The sample $\vec X_1,\dots,\vec X_n$ is generated as previously with~\eqref{eq:ar1} replaced by
\begin{equation}
\tag{EXPAR}
\label{eq:expar}
  X_{ij} = \{ 0.8 - 1.1 \exp ( - 50 X_{i-1,j}^2 ) \} X_{i-1,j} + 0.1 \epsilon_{ij}.
\end{equation}
\end{itemize}

Data under $(\neg H_0) \cap H_{0,m}$, with $H_{0,m}$ as in~\eqref{H0m}, were generated using the procedures described above except that the bivariate random vectors $\vec U_i$, $i \in \{-100,\dots,0,\dots,n\}$ are independent such that $\vec U_i$, $i \in \{-100,\dots,0,\dots,k^\star\}$ are i.i.d.\ from a copula $C_1$ and $\vec U_i$, $i \in \{k^\star+1,\dots,n\}$ are i.i.d.\ from a copula $C_2$, where $C_1 \neq C_2$ and $k^\star = \ip{nt}$ for some $t \in (0,1)$. The resulting samples $\vec X_1,\dots,\vec X_n$ are therefore not samples under $H_{1,c} \cap H_{0,m}$ since the change in the dependence is gradual by~\eqref{eq:ar1} or~\eqref{eq:expar}. The copulas $C_1$ and $C_2$ were taken to be both either bivariate Clayton, Gumbel--Hougaard, Normal or Frank copulas such that $C_1$ has a Kendall's tau of 0.2 and $C_2$ a Kendall's tau of $\tau \in \{0.4,0.6\}$. The parameter $t$ defining $k^\star$ was chosen in $\{0.25, 0.5\}$. 

The dependent multiplier sequences necessary to carry out the tests were generated using the ``moving average approach'' proposed initially in \citet[Section~6.2]{Buh93} and revisited in some detail in \citet[Section~6.1]{BucKoj13}. A standard normal sequence was used for the required initial i.i.d.\ sequence. The kernel function $\kappa$ in that procedure was chosen to be the Parzen kernel defined by $\kappa_{P}(x) = (1 - 6x^2 + 6|x|^3) \1(|x| \leq 1/2) + 2(1-|x|)^3\1(1/2 < |x| \leq 1)$, $x \in \R$, which amounts to choosing the function $\varphi$ in Condition~(M3) as $x \mapsto (\kappa_P \star \kappa_P)(2x) / (\kappa_P \star \kappa_P)(0)$, where `$\star$' denotes the convolution operator. The value of the bandwidth parameter $\ell_n$ defined in Condition~(M2) was chosen using the procedure described in \citet[Section~5]{BucKoj13}. Two choices for the ``combining'' function $\psi$ in that procedure were considered: the median and the maximum. Both choices led to similar rejection rates. The results reported 
in Tables~\ref{H0smmax} and~\ref{H1GHsmmax} below are those obtained with $\psi =$~maximum.

\section{Selected results of the simulation study}
% ================================================
\label{app:tables}

Tables~\ref{H0main} up to~\ref{H1GHsmmax} provide partial results of the large-scale Monte Carlo simulation experiment described in Section~\ref{sec:simulations}. All the tests were carried out at the 5\% level of significance.

\begin{table}%[t!]
\centering
\caption{Percentage of rejection of $H_0$ computed from 1000 random samples of size $n \in \{50, 100, 200\}$ generated under $H_0$, where $C$ is either the $d$-dimensional Clayton (Cl), the Gumbel--Hougaard (GH) or the normal (N) copula whose bivariate margins have a Kendall's tau of $\tau$.} 
\label{H0main}
\begin{tabular}{rrrrrrrrrrrr}
  \hline
  \multicolumn{3}{c}{} & \multicolumn{3}{c}{Cl} & \multicolumn{3}{c}{GH} & \multicolumn{3}{c}{N} \\ \cmidrule(lr){4-6} \cmidrule(lr){7-9} \cmidrule(lr){10-12} $d$ & $n$ & $\tau$ & $\check S_n$ & $\hat S_n$  & $S_n^{R}$ & $\check S_n$ & $\hat S_n$  & $S_n^{R}$ & $\check S_n$ & $\hat S_n$  & $S_n^{R}$  \\ \hline
2 & 50 & 0.00 & 6.2 & 4.0 & 4.6 & 5.4 & 2.9 & 4.5 & 7.3 & 3.4 & 4.8 \\ 
   &  & 0.25 & 6.7 & 6.2 & 5.6 & 5.5 & 3.3 & 5.4 & 4.4 & 3.0 & 6.3 \\ 
   &  & 0.50 & 5.6 & 7.9 & 6.0 & 4.4 & 3.3 & 4.6 & 4.4 & 5.3 & 4.9 \\ 
   &  & 0.75 & 6.0 & 16.6 & 5.5 & 3.2 & 6.7 & 4.3 & 3.6 & 9.1 & 4.9 \\ 
   & 100 & 0.00 & 4.9 & 3.5 & 5.3 & 5.2 & 4.1 & 5.5 & 4.3 & 2.8 & 5.5 \\ 
   &  & 0.25 & 6.1 & 6.6 & 5.0 & 5.0 & 3.3 & 6.2 & 5.3 & 4.0 & 5.5 \\ 
   &  & 0.50 & 4.4 & 9.3 & 5.9 & 3.7 & 2.8 & 5.7 & 3.1 & 3.4 & 5.3 \\ 
   &  & 0.75 & 2.7 & 10.0 & 4.6 & 2.5 & 4.7 & 4.4 & 2.1 & 6.0 & 5.6 \\ 
   & 200 & 0.00 & 4.0 & 3.5 & 5.2 & 4.3 & 4.0 & 5.2 & 5.4 & 4.9 & 4.3 \\ 
   &  & 0.25 & 4.7 & 5.2 & 6.3 & 3.3 & 3.0 & 3.8 & 4.0 & 3.9 & 5.2 \\ 
   &  & 0.50 & 5.1 & 8.5 & 4.9 & 3.2 & 2.3 & 4.5 & 4.0 & 4.7 & 4.8 \\ 
   &  & 0.75 & 2.6 & 9.3 & 5.9 & 1.5 & 3.1 & 5.2 & 1.9 & 4.8 & 5.7 \\ 
  3 & 50 & 0.00 & 4.3 & 1.5 & 3.0 & 4.2 & 2.1 & 3.6 & 5.5 & 2.8 & 3.4 \\ 
   &  & 0.25 & 6.3 & 5.0 & 5.1 & 5.5 & 1.0 & 5.1 & 5.3 & 3.0 & 4.3 \\ 
   &  & 0.50 & 8.2 & 9.1 & 5.9 & 2.7 & 0.9 & 5.7 & 3.0 & 2.2 & 4.6 \\ 
   &  & 0.75 & 2.0 & 2.9 & 6.9 & 0.5 & 0.4 & 6.3 & 1.1 & 1.3 & 4.1 \\ 
   & 100 & 0.00 & 4.5 & 3.4 & 4.5 & 4.5 & 2.8 & 4.6 & 4.5 & 2.7 & 3.9 \\ 
   &  & 0.25 & 5.0 & 5.1 & 5.4 & 4.2 & 2.6 & 4.4 & 5.4 & 3.5 & 4.5 \\ 
   &  & 0.50 & 5.7 & 7.6 & 6.3 & 3.3 & 1.3 & 5.0 & 3.2 & 3.1 & 3.9 \\ 
   &  & 0.75 & 2.5 & 4.9 & 5.0 & 1.0 & 1.0 & 5.2 & 0.8 & 1.6 & 5.5 \\ 
   & 200 & 0.00 & 3.3 & 2.5 & 4.3 & 3.5 & 3.2 & 4.3 & 4.8 & 4.0 & 4.7 \\ 
   &  & 0.25 & 6.6 & 7.1 & 5.5 & 5.0 & 3.3 & 4.5 & 4.8 & 4.1 & 5.0 \\ 
   &  & 0.50 & 6.0 & 9.2 & 4.5 & 3.0 & 2.4 & 5.9 & 4.8 & 4.3 & 4.8 \\ 
   &  & 0.75 & 2.9 & 6.4 & 6.4 & 0.7 & 0.9 & 3.8 & 1.3 & 2.2 & 4.9 \\ 
   \hline
\end{tabular}
\end{table}

\begin{table}%[t!]
\centering
\caption{Percentage of rejection of $H_0$ computed from 1000 samples of size $n \in \{100, 200\}$ generated under $H_0$ as explained in Appendix~\ref{strongmixing}, where $C$ is either the bivariate Clayton (Cl), the Gumbel--Hougaard (GH), the normal (N) or the Frank (F) copula with a Kendall's tau of $\tau$. The columns $\overline{\hat \ell_n^{opt}}$ and std give the mean and the standard deviation of the values of $\ell_n$ used for creating the dependent multiplier sequences.} 
\label{H0smmax}
\begin{tabular}{lrrrrrrrrrr}
  \hline
  \multicolumn{3}{c}{} & \multicolumn{4}{c}{AR1} & \multicolumn{4}{c}{EXPAR} \\ \cmidrule(lr){4-7} \cmidrule(lr){8-11} $C$ & $n$ & $\tau$ & $\overline{\hat \ell_n^{opt}}$ & std & $\check S_n$ & $\hat S_n$ &  $\overline{\hat \ell_n^{opt}}$ & std & $\check S_n$ & $\hat S_n$ \\ \hline
Cl & 100 & 0.00 & 14.2 & 8.5 & 4.2 & 0.7 & 16.7 & 9.5 & 5.1 & 0.6 \\ 
   &  & 0.25 & 14.1 & 8.6 & 5.9 & 1.7 & 16.7 & 10.2 & 5.5 & 2.3 \\ 
   &  & 0.50 & 14.0 & 10.3 & 4.5 & 2.9 & 16.4 & 10.7 & 6.1 & 3.2 \\ 
   &  & 0.75 & 13.3 & 10.0 & 2.7 & 3.3 & 15.5 & 10.7 & 5.1 & 4.7 \\ 
   & 200 & 0.00 & 16.7 & 8.0 & 5.1 & 2.6 & 20.9 & 9.5 & 4.5 & 1.7 \\ 
   &  & 0.25 & 16.0 & 7.3 & 5.1 & 3.0 & 20.5 & 9.7 & 4.1 & 2.0 \\ 
   &  & 0.50 & 15.8 & 7.8 & 2.6 & 2.5 & 19.8 & 9.9 & 3.5 & 2.8 \\ 
   &  & 0.75 & 15.5 & 9.0 & 1.6 & 3.6 & 19.0 & 9.0 & 4.3 & 4.7 \\ 
  GH & 100 & 0.00 & 14.5 & 9.7 & 4.6 & 0.9 & 16.9 & 8.2 & 4.4 & 0.6 \\ 
   &  & 0.25 & 14.1 & 8.7 & 4.9 & 1.5 & 17.1 & 10.2 & 5.0 & 0.6 \\ 
   &  & 0.50 & 14.0 & 9.5 & 3.9 & 1.2 & 15.8 & 9.5 & 4.0 & 0.3 \\ 
   &  & 0.75 & 13.7 & 10.0 & 2.6 & 0.5 & 15.2 & 9.2 & 1.6 & 0.4 \\ 
   & 200 & 0.00 & 16.8 & 7.9 & 4.3 & 1.8 & 21.5 & 10.1 & 3.6 & 1.5 \\ 
   &  & 0.25 & 16.6 & 8.8 & 5.5 & 2.0 & 20.9 & 11.5 & 5.1 & 1.1 \\ 
   &  & 0.50 & 15.9 & 7.4 & 3.7 & 1.6 & 20.1 & 10.9 & 2.9 & 0.7 \\ 
   &  & 0.75 & 15.5 & 8.7 & 1.3 & 0.9 & 18.8 & 9.1 & 1.6 & 0.2 \\ 
  N & 100 & 0.00 & 14.1 & 7.9 & 5.0 & 1.1 & 17.3 & 9.2 & 5.4 & 1.4 \\ 
   &  & 0.25 & 13.5 & 8.1 & 5.9 & 1.4 & 17.3 & 10.8 & 5.0 & 1.1 \\ 
   &  & 0.50 & 13.5 & 9.1 & 3.3 & 1.4 & 16.4 & 9.7 & 4.5 & 1.0 \\ 
   &  & 0.75 & 12.9 & 7.9 & 1.7 & 1.7 & 15.7 & 10.8 & 2.7 & 1.1 \\ 
   & 200 & 0.00 & 16.3 & 6.2 & 5.4 & 1.9 & 20.7 & 8.7 & 3.7 & 1.5 \\ 
   &  & 0.25 & 16.0 & 7.1 & 4.2 & 2.4 & 20.9 & 8.9 & 5.0 & 1.5 \\ 
   &  & 0.50 & 16.1 & 7.8 & 4.2 & 3.2 & 19.8 & 10.4 & 2.9 & 1.8 \\ 
   &  & 0.75 & 15.4 & 7.5 & 0.9 & 1.4 & 19.3 & 10.7 & 0.8 & 0.5 \\ 
  F & 100 & 0.00 & 13.8 & 7.8 & 5.8 & 1.5 & 17.4 & 9.7 & 6.4 & 0.8 \\ 
   &  & 0.25 & 14.2 & 9.3 & 5.5 & 1.3 & 16.6 & 9.7 & 4.8 & 0.5 \\ 
   &  & 0.50 & 13.9 & 9.0 & 3.3 & 1.7 & 16.8 & 11.4 & 3.6 & 0.5 \\ 
   &  & 0.75 & 13.5 & 9.7 & 1.5 & 0.6 & 15.8 & 10.3 & 3.1 & 1.2 \\ 
   & 200 & 0.00 & 16.8 & 7.2 & 4.2 & 2.3 & 20.9 & 9.4 & 4.3 & 1.4 \\ 
   &  & 0.25 & 16.0 & 7.0 & 6.0 & 2.9 & 20.6 & 8.7 & 3.6 & 1.1 \\ 
   &  & 0.50 & 16.1 & 8.2 & 3.1 & 1.5 & 20.3 & 10.4 & 3.2 & 1.2 \\ 
   &  & 0.75 & 15.8 & 8.8 & 0.9 & 0.5 & 19.6 & 9.0 & 1.2 & 0.8 \\ 
   \hline
\end{tabular}
\end{table}

\begin{table}%[t!]
\centering
\caption{Percentage of rejection of $H_0$ computed from 1000 samples of size $n \in \{50, 100, 200\}$ generated under $H_{0,m} \cap H_{1,c}$, where $H_{1,c}$ is defined in~\eqref{H1c}, $k^\star = \ip{nt}$, $C_1$ and $C_2$ are both either bivariate Clayton (Cl), Gumbel--Hougaard (GH) or normal (N) copulas such that $C_1$ has a Kendall's tau of 0.2 and $C_2$ a Kendall's tau of $\tau$.} 
\label{H1onecopd2}
\begin{tabular}{rrrrrrrrrrrr}
  \hline
  \multicolumn{3}{c}{} & \multicolumn{3}{c}{Cl} & \multicolumn{3}{c}{GH} & \multicolumn{3}{c}{N} \\ \cmidrule(lr){4-6} \cmidrule(lr){7-9} \cmidrule(lr){10-12} $n$ & $\tau$ & $t$ & $\check S_n$ & $\hat S_n$  & $S_n^{R}$ & $\check S_n$ & $\hat S_n$  & $S_n^{R}$ & $\check S_n$ & $\hat S_n$  & $S_n^{R}$ \\ \hline
50 & 0.4 & 0.10 & 7.5 & 8.1 & 5.7 & 6.1 & 3.7 & 4.3 & 6.1 & 4.8 & 4.3 \\ 
   &  & 0.25 & 12.1 & 10.6 & 4.0 & 9.4 & 5.1 & 5.0 & 10.4 & 7.6 & 5.5 \\ 
   &  & 0.50 & 18.0 & 16.1 & 6.3 & 12.1 & 7.8 & 4.8 & 12.3 & 8.4 & 5.8 \\ 
   & 0.6 & 0.10 & 14.5 & 17.0 & 5.4 & 11.4 & 7.7 & 6.2 & 11.4 & 9.8 & 7.2 \\ 
   &  & 0.25 & 35.5 & 34.4 & 7.4 & 29.9 & 21.4 & 6.4 & 31.3 & 21.4 & 7.1 \\ 
   &  & 0.50 & 47.3 & 41.6 & 7.0 & 45.3 & 30.3 & 8.9 & 46.0 & 33.9 & 9.1 \\ 
  100 & 0.4 & 0.10 & 7.1 & 8.7 & 6.1 & 6.6 & 5.1 & 5.1 & 5.8 & 5.2 & 5.3 \\ 
   &  & 0.25 & 18.8 & 19.9 & 5.2 & 16.9 & 13.2 & 5.8 & 14.9 & 12.5 & 6.1 \\ 
   &  & 0.50 & 26.5 & 26.4 & 7.3 & 23.8 & 18.8 & 7.3 & 22.6 & 19.1 & 7.9 \\ 
   & 0.6 & 0.10 & 21.5 & 25.1 & 5.8 & 16.7 & 12.0 & 5.1 & 17.5 & 16.9 & 6.1 \\ 
   &  & 0.25 & 65.1 & 66.0 & 6.3 & 61.2 & 51.5 & 7.5 & 62.9 & 54.8 & 9.9 \\ 
   &  & 0.50 & 82.1 & 81.6 & 14.7 & 78.8 & 69.7 & 14.9 & 79.1 & 73.7 & 11.8 \\ 
  200 & 0.4 & 0.10 & 11.1 & 13.8 & 5.8 & 8.3 & 8.1 & 5.5 & 9.5 & 9.8 & 5.0 \\ 
   &  & 0.25 & 30.8 & 33.9 & 5.9 & 27.6 & 24.8 & 6.4 & 29.6 & 28.3 & 6.8 \\ 
   &  & 0.50 & 47.1 & 48.6 & 9.0 & 45.8 & 41.4 & 8.7 & 47.1 & 46.1 & 9.4 \\ 
   & 0.6 & 0.10 & 36.4 & 41.3 & 6.8 & 34.4 & 31.7 & 7.1 & 36.0 & 36.3 & 6.7 \\ 
   &  & 0.25 & 92.6 & 93.2 & 12.3 & 91.4 & 88.9 & 16.7 & 91.3 & 90.2 & 12.0 \\ 
   &  & 0.50 & 98.9 & 99.3 & 22.2 & 98.5 & 98.1 & 22.0 & 99.3 & 99.1 & 21.1 \\ 
   \hline
\end{tabular}
\end{table}

\begin{table}%[t!]
\centering
\caption{Percentage of rejection of $H_0$ computed from 1000 samples of size $n \in \{100, 200\}$ generated under $H_{0,m} \cap H_{1,c}$, where $H_{1,c}$ is defined in~\eqref{H1c}, $k^\star = \ip{nt}$, $C_1$ (resp.\ $C_2$) is a $d$-dimensional Clayton (resp.\ Gumbel--Hougaard) copula whose bivariate margins have a Kendall's tau of $\tau$.} 
\label{H1both}
\begin{tabular}{rrrrrrrrr}
  \hline
  \multicolumn{3}{c}{} & \multicolumn{3}{c}{$d=2$} & \multicolumn{3}{c}{$d=3$} \\ \cmidrule(lr){4-6} \cmidrule(lr){7-9} $n$ & $\tau$ & $t$ & $\check S_n$ & $\hat S_n$  & $S_n^{R}$ & $\check S_n$ & $\hat S_n$  & $S_n^{R}$ \\ \hline
100 & 0.25 & 0.25 & 5.7 & 4.1 & 5.3 & 5.3 & 3.1 & 4.4 \\ 
   &  & 0.50 & 6.2 & 5.7 & 5.6 & 9.1 & 6.0 & 5.8 \\ 
   &  & 0.75 & 6.6 & 6.3 & 3.5 & 5.8 & 5.4 & 4.9 \\ 
   & 0.50 & 0.25 & 5.5 & 5.9 & 5.8 & 4.6 & 2.9 & 5.0 \\ 
   &  & 0.50 & 10.5 & 12.2 & 5.1 & 15.1 & 15.1 & 6.8 \\ 
   &  & 0.75 & 8.3 & 11.9 & 4.4 & 7.7 & 9.9 & 5.3 \\ 
   & 0.75 & 0.25 & 4.0 & 7.6 & 5.1 & 2.5 & 1.9 & 4.3 \\ 
   &  & 0.50 & 12.5 & 19.9 & 6.0 & 9.8 & 13.2 & 5.2 \\ 
   &  & 0.75 & 8.2 & 16.4 & 6.5 & 4.8 & 6.7 & 3.8 \\ 
  200 & 0.25 & 0.25 & 5.8 & 5.3 & 6.1 & 5.8 & 3.7 & 5.8 \\ 
   &  & 0.50 & 8.5 & 8.7 & 5.4 & 9.8 & 9.7 & 6.6 \\ 
   &  & 0.75 & 6.4 & 6.9 & 5.5 & 9.0 & 9.1 & 5.3 \\ 
   & 0.50 & 0.25 & 10.4 & 11.5 & 6.8 & 12.6 & 10.1 & 6.7 \\ 
   &  & 0.50 & 30.9 & 37.5 & 5.3 & 44.0 & 45.8 & 7.1 \\ 
   &  & 0.75 & 16.3 & 23.2 & 6.1 & 20.1 & 27.0 & 5.1 \\ 
   & 0.75 & 0.25 & 11.3 & 18.0 & 4.9 & 15.6 & 16.7 & 7.3 \\ 
   &  & 0.50 & 43.4 & 54.4 & 6.2 & 58.9 & 63.1 & 5.2 \\ 
   &  & 0.75 & 21.3 & 36.4 & 4.6 & 23.6 & 36.0 & 6.8 \\ 
   \hline
\end{tabular}
\end{table}

\begin{table}%[t!]
\centering
\caption{Rejection percentage of $H_0$ computed from 1000 samples of size $n \in \{50, 100, 200\}$ such that the $\ip{nt_1}$ first observations of each sample are from a $d$-variate c.d.f.\ with normal copula and $N(0,1)$ margins (that is, from a multivariate standard normal c.d.f.), and the $n - \ip{nt_1}$ last observations are from a $d$-variate c.d.f.\ with normal copula whose first margin is the $N(\mu,1)$ and whose $d-1$ remaining margins are the $N(0,1)$. The bivariate margins of the normal copula have a Kendall's tau of~$\tau$.} 
\label{H0break}
\begin{tabular}{rrrrrrrrrrrrrrr}
  \hline
  \multicolumn{3}{r}{$(\mu,t_1)=$} & \multicolumn{3}{c}{$(0.5,0.25)$} & \multicolumn{3}{c}{$(0.5,0.5)$} & \multicolumn{3}{c}{$(2,0.25)$} & \multicolumn{3}{c}{$(2,0.5)$} \\ \cmidrule(lr){4-6} \cmidrule(lr){6-9} \cmidrule(lr){10-12} \cmidrule(lr){13-15} $d$ & $n$ & $\tau$ & $\check S_n$ & $\hat S_n$  & $S_n^{R}$ & $\check S_n$ & $\hat S_n$  & $S_n^{R}$ & $\check S_n$ & $\hat S_n$  & $S_n^{R}$ & $\check S_n$ & $\hat S_n$ & $S_n^{R}$ \\ \hline
2 & 50 & 0.00 & 6.7 & 3.8 & 9.0 & 6.5 & 3.1 & 17.7 & 4.6 & 2.4 & 70.1 & 5.1 & 2.6 & 98.5 \\ 
   &  & 0.25 & 5.6 & 2.8 & 8.1 & 4.8 & 3.6 & 15.0 & 6.0 & 3.1 & 57.7 & 4.7 & 1.6 & 98.5 \\ 
   &  & 0.50 & 4.4 & 4.1 & 9.6 & 3.4 & 2.9 & 13.0 & 18.4 & 5.1 & 39.6 & 5.7 & 0.9 & 99.3 \\ 
   & 100 & 0.00 & 5.3 & 4.0 & 19.0 & 6.3 & 5.2 & 30.4 & 6.1 & 4.1 & 99.0 & 5.5 & 2.5 & 100.0 \\ 
   &  & 0.25 & 5.0 & 4.4 & 13.9 & 3.7 & 2.3 & 24.9 & 8.6 & 4.7 & 97.5 & 4.4 & 2.2 & 100.0 \\ 
   &  & 0.50 & 4.0 & 3.8 & 11.5 & 3.1 & 3.1 & 22.6 & 29.5 & 13.1 & 91.6 & 12.4 & 2.0 & 100.0 \\ 
   & 200 & 0.00 & 5.5 & 5.2 & 34.1 & 3.9 & 3.4 & 61.9 & 4.3 & 3.3 & 100.0 & 5.6 & 4.1 & 100.0 \\ 
   &  & 0.25 & 3.9 & 3.4 & 27.6 & 4.1 & 3.5 & 51.3 & 13.6 & 8.9 & 100.0 & 8.1 & 4.0 & 100.0 \\ 
   &  & 0.50 & 3.4 & 3.9 & 18.9 & 2.8 & 2.9 & 43.1 & 57.8 & 39.4 & 100.0 & 34.8 & 8.3 & 100.0 \\ 
  3 & 50 & 0.00 & 4.9 & 1.6 & 5.0 & 4.5 & 1.7 & 10.0 & 4.8 & 1.9 & 36.5 & 6.0 & 2.3 & 79.6 \\ 
   &  & 0.25 & 5.0 & 2.7 & 6.9 & 5.0 & 2.6 & 9.9 & 6.9 & 2.6 & 24.8 & 5.1 & 2.6 & 87.4 \\ 
   &  & 0.50 & 3.7 & 2.6 & 5.5 & 3.3 & 1.3 & 8.9 & 9.4 & 3.5 & 17.4 & 3.0 & 0.7 & 94.2 \\ 
   & 100 & 0.00 & 4.5 & 2.3 & 11.3 & 4.2 & 2.5 & 18.6 & 4.8 & 3.0 & 87.5 & 3.6 & 2.0 & 99.6 \\ 
   &  & 0.25 & 4.9 & 3.1 & 10.7 & 5.1 & 3.5 & 14.8 & 6.7 & 3.8 & 67.9 & 5.9 & 4.2 & 99.9 \\ 
   &  & 0.50 & 2.8 & 2.0 & 7.3 & 3.0 & 2.4 & 13.7 & 16.9 & 8.6 & 60.6 & 6.0 & 1.3 & 100.0 \\ 
   & 200 & 0.00 & 3.0 & 2.4 & 20.1 & 4.5 & 4.0 & 37.3 & 5.3 & 3.3 & 100.0 & 4.3 & 3.2 & 100.0 \\ 
   &  & 0.25 & 4.8 & 4.1 & 15.3 & 5.6 & 4.4 & 30.6 & 11.9 & 8.1 & 99.2 & 7.4 & 5.2 & 100.0 \\ 
   &  & 0.50 & 4.9 & 3.7 & 11.8 & 3.8 & 3.2 & 24.9 & 41.0 & 30.9 & 99.2 & 23.9 & 7.8 & 100.0 \\ 
   \hline
\end{tabular}
\end{table}

\begin{table}%[t!]
\centering
\caption{Percentage of rejection of $H_0$ computed from 1000 samples of size $n \in \{100, 200\}$ generated under $\neg H_0$ as explained in the second paragraph of Appendix~\ref{strongmixing}, where $C_1$ and $C_2$ are both bivariate Gumbel--Hougaard copulas such that $C_1$ has a Kendall's tau of 0.2 and $C_2$ a Kendall's tau of $\tau$. The columns $\overline{\hat \ell_n^{opt}}$ and std give the mean and the standard deviation of the values of $\ell_n$ used for creating the dependent multiplier sequences.} 
\label{H1GHsmmax}
\begin{tabular}{rrrrrrrrrrr}
  \hline
  \multicolumn{3}{c}{} & \multicolumn{4}{c}{AR1} & \multicolumn{4}{c}{EXPAR} \\ \cmidrule(lr){4-7} \cmidrule(lr){8-11} $n$ & $t$ & $\tau$ & $\overline{\hat \ell_n^{opt}}$ & std & $\check S_n$ & $\hat S_n$ &  $\overline{\hat \ell_n^{opt}}$ & std & $\check S_n$ & $\hat S_n$ \\ \hline
100 & 0.25 & 0.4 & 14.2 & 10.1 & 13.8 & 3.8 & 16.9 & 10.5 & 10.6 & 1.5 \\ 
   &  & 0.6 & 14.1 & 9.1 & 39.9 & 10.6 & 15.8 & 9.5 & 32.6 & 8.5 \\ 
   & 0.50 & 0.4 & 14.3 & 10.1 & 18.0 & 5.3 & 17.1 & 10.3 & 13.5 & 2.5 \\ 
   &  & 0.6 & 13.9 & 8.6 & 57.2 & 25.2 & 16.5 & 9.3 & 54.8 & 17.1 \\ 
  200 & 0.25 & 0.4 & 16.5 & 8.6 & 16.4 & 8.6 & 20.6 & 9.8 & 14.6 & 4.5 \\ 
   &  & 0.6 & 16.4 & 8.0 & 71.2 & 46.0 & 19.6 & 9.0 & 63.1 & 32.5 \\ 
   & 0.50 & 0.4 & 16.6 & 7.4 & 31.9 & 19.0 & 20.7 & 8.8 & 25.8 & 12.4 \\ 
   &  & 0.6 & 16.3 & 7.1 & 89.8 & 75.5 & 20.6 & 9.1 & 80.9 & 62.5 \\ 
   \hline
\end{tabular}
\end{table}

\end{document}